\documentclass[10pt, a4paper]{article}

\usepackage{amssymb}
\usepackage{amsthm}
\usepackage{mathtools, amsmath, amsfonts, bbm}
\usepackage[english]{babel}
\usepackage{algorithmic}
\usepackage[latin1]{inputenc}
\usepackage{pifont, fontenc, tabularx}
\usepackage{multirow,url, multicol}
\usepackage[usenames,dvipsnames]{xcolor}
\usepackage{tikz}
\usepackage{dsfont, dutchcal}

\newcommand{\R}{\mathbb{R}}
\newcommand{\N}{\mathbb{N}}

\newcommand{\Id}{{\mathds{1}}}
\newcommand{\tr}{\mathrm{tr}} 
\renewcommand{\d}{{\,\mathrm{d}}} 
\newcommand{\dist}{\mathrm{dist}}
\newcommand{\energy}{\mathcal{W}}
\newcommand{\pathlength}{\mathcal{L}}
\newcommand{\pathenergy}{\mathcal{E}}
\newcommand{\Pathenergy}{\mathbf{E}}
\newcommand{\average}{\mathcal{Av}} 
\newcommand{\Average}{\mathbf{Av}} 
\newcommand{\interpolation}{\mathcal{I}} 
\newcommand{\Interpolation}{\mathbf{I}}
\newcommand{\metric}{g}
\newcommand{\ptransport}{{\mathrm{P}}} \newcommand{\Ptransport}{{\mathbf{P}}} \newcommand{\shell}{s}
\newcommand{\shellspace}{\mathcal{S}}
\newcommand\calB{{\mathcal{B}}}
\newcommand\calH{{\mathcal{H}}}
\newcommand\CS{{\mathcal{CS}}}
\newcommand{\discB}[1]{{\mathbf{B}^{#1}}}
\newcommand{\discH}[1]{{\mathbf{H}^{#1}}}
\newcommand{\discCS}[1]{{\mathbf{CS}^{#1}_{\kappa}}}
\newcommand{\control}{d}
\newcommand*\rfrac[2]{{}^{#1}\!/_{#2}}
\newcommand{\discExp}{{\mathbf{Exp}}} \newcommand{\discLog}{{\mathbf{Log}}} \newcommand{\mem}{{\mbox{{\tiny mem}}}}
\newcommand{\bend}{{\mbox{{\tiny bend}}}}

\DeclarePairedDelimiter\floor{\lfloor}{\rfloor}

\DeclareMathOperator{\Hess}{{Hess}}
\DeclareMathOperator*{\argmin}{argmin}
\DeclareMathOperator{\firstFF}{{I}}
\DeclareMathOperator{\scndFF}{{II}}

\definecolor{red}{HTML}{E41A1C}
\definecolor{blue}{HTML}{377EB8}
\definecolor{green}{HTML}{4DAF4A}
\definecolor{purple}{HTML}{984EA3}

\providecommand{\keywords}[1]{\noindent\textbf{Keywords} #1}
\providecommand{\msc}[1]{\noindent\textbf{MSC(2010)} #1}

\title{Smooth Interpolation of Key Frames in a Riemannian Shell Space}
\author{Pascal Huber, Ricardo Perl and Martin Rumpf\footnote{Institut f\"ur Numerische Simulation,
  Universit\"at Bonn, 53115 Bonn (Germany)}}

\begin{document}

\maketitle

\begin{abstract}
  \noindent Splines and subdivision curves are flexible tools in the design and
  manipulation of curves in Euclidean space.  In this paper we study
  generalizations of interpolating splines and subdivision schemes to the
  Riemannian manifold of shell surfaces in which the associated metric measures
  both bending and membrane distortion. The shells under consideration are
  assumed to be represented by Loop subdivision surfaces.  This enables the
  animation of shells via the smooth interpolation of a given set of key frame
  control meshes.  Using a variational time discretization of geodesics
  efficient numerical implementations can be derived. These are based on a
  discrete geodesic interpolation, discrete geometric logarithm, discrete
  exponential map, and discrete parallel transport.  With these building blocks
  at hand discrete Riemannian cardinal splines and three different types of
  discrete, interpolatory subdivision schemes are defined.  Numerical results
  for two different subdivision shell models underline the potential of this
  approach in key frame animation.
\end{abstract}
\bigskip


\medskip
\msc{68U07, 65D05, 65D17, 49M25}

\section{Introduction}
\label{sec:introduction}

In this paper we investigate classical interpolation and approximation tools for
curves in the context of Riemannian manifolds and specifically a shape manifold
of subdivision shell surfaces.  To this end we pick up the \emph{time-discrete
  geodesic calculus} developed by Rumpf and Wirth in~\cite{RuWi12b} and applied
to shell space by Heeren et al.\ in~\cite{HeRuWa12,HeRuSc14}, and ask for
effective and efficient generalizations of the \emph{de Casteljau algorithm} for
\emph{B\'ezier curves}, \emph{cardinal splines}, and several \emph{interpolatory
  subdivision schemes}.

Recently, Riemannian calculus on shape spaces and in particular on the space of
surfaces has attracted a lot of attention.  Kilian et al.~\cite{KiMiPo07}
studied geodesics in the space of triangular surfaces to interpolate between two
given poses.  The underlying metric is derived from the in-plane membrane
distortion.  Since this pioneering paper a variety of other Riemannian metrics
on the space of surfaces have been investigated
\cite{LiShDi10,KuKlDi10,KuKlGo12,BaHaMi11,JiZeLu09,BaBr11,BrBrKi09,AdOvWa08}.
In~\cite{HeRuWa12,HeRuSc14} a metric was proposed that takes the full elastic
responses including bending distortion into account.  We pick up this metric and
the associated variational discretization, which are briefly revisited below.
Also based on this metric, Brandt et al.~\cite{BrTyHi16} proposed an accelerated
scheme for the computation of geodesic paths in shell space. To this end they
applied dimension reduction with respect to the surface models.  Winkler et
al.~\cite{WiDrAl10} and~Fr{\"o}hlich and Botsch~\cite{FrBo11} considered a
representation of triangular meshes via a vector of edge lengths and dihedral
angles and applied a back projection onto the space of triangular surfaces.
This approach has been used by Heeren et al in~\cite{HeRuSc16} to simplify the
computation of Riemannian splines in the space of images.  For the smooth
interpolation of given key frames and independent from the Riemannian concept,
spacetime constraints were used by Witkin and Kass~\cite{WiKa88} to compute
optimal trajectories. Smooth interpolating paths in high dimensional shape
spaces have been investigated via a spacetime constraint approach as
generalizations of traditional cubic B-splines by Kass and Anderson named
\emph{wiggly splines}~\cite{KaAn08}.  These splines lack the theoretical
foundation of an underlying Riemannian shape manifold.  Another alternative
approach is the animation approach by Schulz et al.~\cite{ScTySe14,ScTySe15},
which combines the wiggly spline method, suitable vibration modes and a new
warping scheme.  A discrete Riemannian approach for the smoothing of curves in
a shape space of triangular meshes has been proposed by Brandt et
al.~\cite{BrTyHi16}. They apply a gradient descent scheme for a discrete path
energy on shell space picking up the the concept of the discrete path energy
introduced in~\cite{HeRuWa12}.

In Euclidean space cubic splines minimize the total squared second derivative.
Analogously, \emph{Riemannian cubic splines} were introduced by Noakes et
al.~\cite{NoHePa89} in a variational setting on Riemannian manifolds as smooth
curves that are stationary paths of the integrated squared covariant derivative
of the velocity.  Subsequently, Camarinha et al.~\cite{CaSiCr01} proved a local
optimality condition of this functional at a critical point and Giambo and
Giannoni~\cite{GiGi02} established a global existence results. A typical
application of this kind of higher-order interpolation is for instance path
planning as it occurs naturally in aerospace and manufacturing industries.  The
associated Euler-Lagrange equation is a fourth-order differential equation,
first derived in Noakes et al.~\cite{NoHePa89} and then in Crouch and Silva
Leite~\cite{CrSi95} in the context of dynamic interpolation. Crouch and Silva
Leite also considered the interpolation problem for multiple points on a
manifold.  More recently, Trouv\'{e} and Vialard~\cite{TrVi12} developed a
spline interpolation method on Riemannian manifolds and applied it to
time-indexed sequences of 2D or 3D shapes where they focused on the finite
dimensional case of landmarks.  They introduced a control variable $u$ on the
Hamiltonian equations of the geodesics.  Hinkle et al.~\cite{HiMuFl12}
introduced a family of higher order Riemannian polynomials to perform polynomial
regression on Riemannian manifolds. They apply their approach to low dimensional
embedded manifolds (e.g.\ the $d$-dimensional sphere), to the Lie group $SO(3)$
as well as shape spaces of 2D image data represented by landmark positions.
In~\cite{HeRuSc16} a variational discretization of splines on shell space was
presented. Thereby, discrete splines are defined as minimizers of the time
discrete spline energy coupled with a set of variational constraints. Compared
to the approach considered here, this method leads to a globally coupled system
of all shells along a discrete curve in shell space, which is computationally
challenging.

Different from these methods based on the \emph{intrinsic} formulation using the
covariant derivative on curves, \emph{extrinsic variational formulations} have
been studied, which minimize curve energies in ambient space where the
restriction of the curve to the manifold is realized as a
constraint. Wallner~\cite{Wa04} showed existence of minimizers in this setup for
finite dimensional manifolds, and Pottmann and Hofer~\cite{PoHo05} proved that
these minimizers are $C^2$.  In~\cite{HoPo04} they already provided a method for
the computation of splines on parametric surfaces, level sets, triangle meshes
and point set surfaces.  Algorithmically, they alternately compute minimizers of
the objective functional in the tangent plane and project back to the manifold.
Our approach can not be formulated extrinsically, because the metric on shell
space does not arise from an ambient Euclidean metric.

The concept of B\'ezier curves can easily be transferred to Riemannian manifolds
and in particular to shape spaces.  To this end the de Casteljau algorithm with
linear interpolation has to be replaced by geodesic
interpolation~\cite{PoNo07,GoSaAb14,AbGoSt16}, compare also Effland et
al.~\cite{EfRuSi14} for a corresponding tool on the shape space of images and
Brandt et al.~\cite{BrTyHi16} for this method on the space of triangular shell
surfaces using the variational time discretization~\cite{HeRuWa12}.  Below we
will discuss the de Casteljau algorithm on the space of subdivision surfaces to
prepare the discussion of Hermite interpolation and cardinal splines.

In addition to variational formulations there are various approaches for
interpolatory \emph{subdivision schemes} on manifolds.  These methods exploit
the fact that subdivision schemes for curves in linear spaces are mostly based
on repeated local averages. Dyn~\cite{Dy02} proposed a Riemannian extension,
where affine averaging is replaced by geodesic averaging.  Wallner and
Dyn~\cite{WaDy05} showed that the Riemannian extension of cubic subdivision
yields $C^1$ curves.  Their analysis is based on a comparison of subdivision
schemes in Euclidean space---where $C^1$ smoothness results are obtained by
studying the convergence of the symbol associated with the linear subdivision
rule---and the corresponding scheme on a Riemannian manifold.  Under a suitable
proximity condition $C^1$ smoothness can be established.  In~\cite{Wa06b} this
technique is generalized to show higher smoothness and specifically for Lie
groups $C^2$ smoothness of Riemannian subdivision schemes.  More general
proximity conditions which imply higher smoothness of the Riemannian counterpart
of an linear subdivision scheme are analysed by Grohs in~\cite{Gr10}.  Dyn et
al.\ show in~\cite{DyGr10} that the approximation order of nonlinear univariate
schemes derived as perturbations of linear schemes is directly linked to the
smoothness of the nonlinear scheme.  In~\cite{DuXi13} Duchamp et al.\ studied
higher order smoothness of manifold valued subdivision curves using a retraction
map in the construction of the subdivision scheme.  This retraction map is
required to be a third-order approximation of the exponential map.  Recently,
Wallner~\cite{Wa14} proved convergence of the linear four-point scheme and other
univariate interpolatory schemes in Riemannian manifolds.  To this end he
studies Riemannian edge length contractivity of the schemes on a manifold
combined with a multiresolutional analysis.  In our case, concerning the
consistency, there is not only the transfer of linear subdivision schemes to
curved spaces. But there is also a second source of proximity errors, which is
the approximation of the local squared Riemannian distance by a functional which
is cheaper to compute. At the moment the convergence of the schemes considered
here for fixed time step size is open. We refer to Section~\ref{sec:discussion}
for some comments in this direction.

We consider here Loop's subdivision also for the representation of our shell
surfaces. They are subdivision limit surfaces for given control meshes.
Subdivision schemes for the modeling of geometries are widespread in geometry
processing and computer graphics.  For a comprehensive introduction to
subdivision methods in general we refer the reader to~\cite{PeRe08}
and~\cite{Ca12}.  Among the most popular subdivision schemes are the
Catmull-Clark~\cite{CaCl78} and Doo-Sabin~\cite{DoSa78} schemes on quadrilateral
meshes, and Loop's scheme on triangular meshes~\cite{Lo94}.  The Loop
subdivision basis functions have been extensively studied, in particular
regarding their smoothness~\cite{Re95}, curvature integrability~\cite{ReSc01},
(local) linear independence~\cite{PeWu06a, ZoJKo14}, approximation
power~\cite{Ar01}, and robust evaluation around extraordinary
vertices~\cite{St99}.  With respect to the use of subdivision methods in
animation see~\cite{DeKaTr98, ThWaSt06}.  Nowadays, subdivision finite elements
are also extensively used in engineering~\cite{CiOrSc00, CiLo11, CiScAn02,
  CiOr01, GrTuSt02, GrTu04}.

\paragraph{Our contribution} In this paper we use recent progress on the theory
of splines on Riemannian manifolds and take into account the Riemannian model of
the space of shells proposed by Heeren et al.~\cite{HeRuWa12}.
We show that various Euclidean concepts of interpolating curves can easily
be transferred to the Riemannian shell space. To this end we slightly
generalize the discrete geodesic calculus developed in~\cite{HeRuWa12}
and~\cite{HeRuSc14} which allows for an elegant and robust implementation of a
versatile toolbox for key frame interpolation. In detail our contributions
are 
\begin{itemize}
\item[-] the generalization of discrete geodesics and the discrete exponential
  map introduced in~\cite{HeRuWa12} and~\cite{HeRuSc14} to points in time which
  are not multiples of a given time step size as an essential ingredient for the
  construction of smooth interpolating curves,
\item[-] a discrete Riemannian Catmull-Rom interpolation based on discrete
  B{\'e}zier curves and discrete parallel transport,
\item[-] different discrete Riemannian interpolatory subdivision schemes,
\item[-] and the conforming implementation of the geodesic calculus on the space of
  shells represented by subdivision surfaces including the new interpolation
  curves. 
\end{itemize}
\emph{The organization of the paper is as follows.}  In
Section~\ref{sec:vari-time-discr} we recall some basic facts of Riemannian
calculus and discuss the variational time discretization and the set up in the
case of the space of shells represented by Loop subdivision surfaces.
Section~\ref{sec:bezier-curves-shell} is devoted to the review of discrete
B{\'e}zier curves in shell space, whereas in Section \ref{sec:time-discr-card} we
derive discrete cardinal splines. Then in Section~\ref{sec:interp-subd-curv}
different discrete interpolatory subdivision schemes are presented. Finally, in
Section~\ref{sec:discussion} we draw conclusions.

\section{Variational time discretization of geodesic calculus in shell space}
\label{sec:vari-time-discr}

In this section we review the \emph{variational time discretization of geodesic
  calculus} developed in~\cite{RuWi12b} by Rumpf and Wirth and applied to a
space of discrete shells by Heeren et al.\ in~\cite{HeRuWa12,HeRuSc14}.  This
calculus includes in particular the notions of \emph{discrete geodesics},
\emph{discrete logarithm}, \emph{discrete exponential map}, and \emph{discrete
  parallel transport}.  For results on existence and uniqueness of the
associated geometric operators and the convergence to their continuous
counterparts on finite dimensional Riemannian manifolds and certain infinite
dimensional Hilbert manifolds we refer to~\cite{RuWi12b}.

\paragraph{Geodesic calculus on Riemannian manifolds} 
Let us denote by $(\shellspace,\,\metric)$ a smooth, complete \emph{Riemannian
  manifold} with metric $\metric$.  In our context \(\shellspace\) is a space of
discrete shells and more specifically the space of shells represented by Loop
subdivision surfaces.  For a smooth path $( \shell(t) )_{t \in [0,1]}$, the
\emph{path length} is defined by
\begin{align}
  \pathlength[ ( \shell(t) )_{t \in [0,1]} ] = 
  \int_0^1 \sqrt{\metric_{\shell(t)} ( \dot \shell(t), \dot \shell(t) )} \, dt, \label{eq:pathLength}
\end{align}
where $\dot \shell(t) \in T_{\shell(t)}\shellspace$ denotes the velocity at time
$t$.  Given two shells $\shell_A, \shell_B \in \shellspace$, the minimizing path
$( \shell(t) )_{t \in [0,1]}$ of \eqref{eq:pathLength} is called shortest
\emph{geodesic} and the minimal path length defines the Riemannian distance
$\dist(\shell_A,\shell_B)$ between $\shell_A$ and $\shell_B$.  Geodesic paths
also minimize the \emph{path energy}
\begin{align}
  \pathenergy[ ( \shell(t) )_{t \in [0,1]} ] = \int_0^1 \metric_{\shell(t)} ( \dot \shell(t), \dot \shell(t) ) \, dt, \label{eq:pathEnergy}
\end{align}
if in addition the speed
$\sqrt{\metric_{\shell(t)} ( \dot \shell(t), \dot \shell(t) )}$ is constant. For
small Riemannian distances shortest geodesics are unique.  Then the
\emph{geodesic averaging} for $0 \leq t\leq 1$ is given by
\begin{align}
  \label{eq:1}
  \average( \shell_A, \shell_B, t)  = \shell(t) 
\end{align}
where $\shell(t)$ is assumed to be the shortest geodesic connecting $\shell_A$
and $\shell_B$.  The \emph{exponential map} is defined by
\begin{align*}
  \exp_{\shell_A}(v) = \shell(1)\,, 
\end{align*}
where $t \mapsto \shell(t)$ is the solution of the Euler-Lagrange equation
$\nabla_{\dot \shell(t)} \dot \shell(t) = 0$ associated with the path energy for
given $\shell(0) = \shell_A$ and
$\dot \shell(0) = v \in T_{\shell_A}\shellspace$. Here $\nabla_{\dot \shell(t)}$
denotes the covariant derivative along the geodesic.  For sufficiently small $r$
the exponential map $\exp_{\shell_A}$ is a bijection from
$B_r(0)\to \exp_{\shell_A}(B_r(0))$.  For $\shell_B \in \exp_{\shell_A}(B_r(0))$
the inverse of $\exp_{\shell_A}$ defines the \emph{logarithm map}
\begin{align*}
  \log_{\shell_A}(\shell_B) = v,
\end{align*}
where $v$ is the initial velocity $\dot \shell(0)$ of the unique shortest
geodesic connecting $\shell_A$ with $\shell_B$.  Finally, the \emph{parallel
  transport} of a given vector $w$ along a path $(\shell(t))_{t\in[0,1]}$ (not
necessarily geodesic) is defined as
\begin{align*}
  \ptransport_{(\shell(t))_{t\in[0,1]}}(w) = w(1)\,,
\end{align*}
where $t\to w(t)$ solves $\nabla_{\dot \shell(t)} w(t) = 0$ with initial data
$w(0) = w$.  For a detailed discussion we refer to~\cite{DoCarmo1992}.

\paragraph{Subdivision shell space}
For the manifold $\shellspace$ of discrete shells investigated in this paper an
element of this space is a Loop subdivision surface described by a triangular
control mesh $\shell$.  Thus for fixed topology of the control mesh and a point
to point correspondence of control meshes in $\shellspace$ a discrete shell can
be identified with a vector in $\R^{3M}$, where $M$ is the number of vertices of
the control mesh.  Physically, a subdivision surface with control mesh $\shell$
can be considered as the mid-surface of a thin, curved three-dimensional
volumetric object whose thickness $\delta$ is relatively small compared to the
area of the mid-surface.  The associated metric reflects viscous dissipation in
this thin surface layer when undergoing a deformation.  To derive the metric we
apply Rayleigh's paradigm (cf.~\cite{HeRuWa12}), i.e.\ we start with a model for
the elastic deformation of the thin layer and replace then strain by strain
rates for a second order approximation of the energy.

To this end, we consider at first an elastic deformation $\phi$ of the
homogeneous, isotropic material layer and approximate the energy by an elastic
energy associated to the deformation of the midsurface split into a membrane and
bending contribution; i.e.\ we define
\begin{align*}
  Q_{\mem}[\phi]&=B^{\mem}_{\phi}-B^{\mem}_{\Id}\quad\text{and}\quad
                  Q_{\bend}[\phi]=B^{\bend}_{\phi}-B^{\bend}_{\Id}
\end{align*}
with the symmetric operators $B^{\mem}$ and $B^{\bend}$ given by
\begin{align*}
  \firstFF_{\shell}(B^{\mem}_{\phi}v,w)&=\firstFF_{\phi(\shell)}(\d\phi(v),\d\phi(w))\,,\\
  \firstFF_{\shell}(B^{\bend}_{\phi}v,w)&=\scndFF_{\phi(\shell)}(\d\phi(v),\d\phi(w))\, .
\end{align*}
Here, $\firstFF_\shell$, $\scndFF_\shell$ denote the first and second
fundamental form with tangent vectors $v$ and $w$.  In fact, $B^{\mem}_{\phi}$
is the geometric (tangential) Cauchy-Green strain tensor for the deformation
$\phi$ on the surface $\shell$ and $Q_{\bend}[\phi]$ is the relative shape
operator quantifying properly the difference in curvature between the deformed
and the undeformed configuration.  The total elastic deformation energy is then
given by
\begin{align}
  \label{eq:totalEnergy}
  \energy_\shell^D[\phi]=\delta \int_\shell W_{\mem}(Q_{\mem}[\phi])\,\d x
  + \delta^3\int_\shell W_{\bend}(Q_{\bend}[\phi])\,\d x \ , 
\end{align}
where the non-negative energy densities $W_{\mem}$ and $ W_{\bend}$ act on the
symmetric linear rank two operators \(Q_{\mem}\) and \(Q_{\bend}\),
respectively.  For $W_{\mem}$ and $W_{\bend}$ we require that (i) $W(0) = 0$,
(ii) $DW = 0$ at the zero matrix, and (iii) $D^2 W$ is positive definite at the
zero matrix. These requirements ensure that for a shell in a stress free
configuration the deformation identity $\Id$ is a minimizer of $\energy$ and
thus (i) $\energy_\shell^D[\Id] =0$ and (ii) $d \energy_\shell^D[\Id] =
0$. Additionally, we assume (iii) that the energy is strictly convex (modulo
rigid body motions) in a neighborhood of a minimizer. These conditions capture
most thin elastic materials~\cite{Ci00}.  Under these assumptions the
application of Rayleigh's paradigm indeed leads to a Riemannian metric as
pointed out by Heeren et al. \cite{HeRuSc14}.
They showed the following result :\\[1ex]
{\it For $v$ a tangent vector field to a shell $\shell$ in the space of smooth
  shells $\shellspace$,  $\Hess(\energy_\shell^D)(v,v)=0$ if and only if $v$ induces an
  infinitesimal rigid motion. Consequently,
  $g_\shell(v,w) = \frac12 \Hess(\energy_\shell^D)(v,w)$ is indeed a Riemannian
  metric on the space of
  smooth shells modulo infinitesimal rigid body motions.}\\[1ex]
For the proof we refer to~\cite{HeRuSc14}.  In the application, we consider the
membrane energy density as suggested in~\cite{HeRuWa12}:
\begin{align}
  \label{eq:membraneDensity}
  W_{\mem}(Q) &= \frac{\mu}{2}\tr Q + \frac{\lambda}{4} \det Q - \frac{2\mu+\lambda}{4} \log \det Q - \mu - \frac{\lambda}{4}.
\end{align}
Here, $\lambda$ and $\mu$ are the Lam\'e constants of linearized elasticity and
$\tr $ and $\det $ denote the trace and the determinant, respectively, where
$\det Q$ describes area distortion, while $\tr Q$ measures length distortion.
The function $W_{\mem}[Q]$ is \emph{rigid body motion invariant} and the
identity as deformation is the minimizer.  Furthermore, the Hessian of the
resulting path energy coincides with the quadratic form of linearized
elasticity.  The $\log \det Q$ term penalizes material compression.  For the
bending energy density we select the squared Frobenius norm of the embedded
relative Weingarten map
\begin{align}
  \label{eq:bendingDensity}
  W_\bend(Q)  = \tr( Q^T Q )\,.
\end{align}

\paragraph{Time-discrete geodesic calculus}
To derive a variational discretization of the geodesic calculus we consider an
approximation of the squared Riemannian distance $\dist^2$ by a smooth
functional $\energy: \shellspace \times \shellspace \to \R$, where we assume
$\shell, \tilde \shell \in \shellspace$:
\begin{align}
  \dist^2(\shell,\tilde \shell) = \energy[\shell,\tilde \shell] + \dist^3(\shell,\tilde \shell). \label{eq:approxSqrDist}
\end{align}
For a detailed discussion of this approximation we refer to~\cite{HeRuWa12}.
For $\metric_\shell = \frac{1}{2} \energy_{,22}[\shell,\shell]$ this assumption
on $\energy$ holds by Taylor expansion.  In the case of the subdivision shell
space introduced above,
$\energy[\shell,\tilde \shell] := \energy_\shell^D[\tilde \shell- \shell]$
fulfills the assumption with $\tilde \shell- \shell$ representing the
deformation of $\shell$ to $\tilde \shell$ as a vector-valued subdivision
function on $\shell$. Using the Cauchy-Schwarz inequality one obtains that
$\pathenergy[(\shell(t))_{t\in [0,1]}] \geq \sum_{k=1}^{K}
\dist(\shell(\frac{k-1}{K}),\shell(\frac{k}{K}))^2$, where equality holds only
if $(\shell(t))_{t\in [0,1]}$ is already a shortest geodesic.
This motivates the definition of the \emph{discrete path energy}
\begin{align}
  \Pathenergy[ (\shell_0, \ldots, \shell_K) ] = K \sum_{k=1}^{K} \energy[ \shell_{k-1}, \shell_{k} ] \label{eq:discPathEnergy}
\end{align}
of a \emph{discrete path} $(\shell_0, \ldots, \shell_K)$.  For given input
shells $\shell_A$ and $\shell_B$ in $\shellspace$ we call the discrete path
$(\shell_0, \ldots, \shell_K)$ a \emph{discrete shortest geodesic}, if
$(\shell_0=\shell_A, \ldots, \shell_K=\shell_B)$ is a minimizer of the discrete
path energy \eqref{eq:discPathEnergy}.  The corresponding system of
Euler--Lagrange equations is given by
\begin{align}
  \label{eq:EL}
  \energy_{,2}[\shell_{k-1},\shell_{k}] + \energy_{,1}[\shell_{k},\shell_{k+1}] =0
\end{align}
for $k=1,\ldots, K-1$, where $\energy_{,i}$ denotes the variation with respect
to the $i$th argument.
Then the \emph{discrete geodesic averaging} for $0 \leq k \leq K$ is defined by
\begin{align}
  \label{eq:1}
  \Average^{K}( \shell_A, \shell_B, k/K) = \shell_k 
\end{align}
if $(\shell_0, \ldots, \shell_k,\ldots, \shell_K)$ is the shortest discrete
geodesic connecting $\shell_0=\shell_A$, $\shell_K=\shell_B$.
Given a continuous geodesic $( \shell(t) )_{t \in [0,1]}$ with
$\shell(0) = \shell_A$ and $\shell(1) = \shell_B$ and a discrete geodesic
$(\shell_0, \ldots, \shell_K)$ with $\shell_0 = \shell_A$ and
$\shell_K = \shell_B$, we may view $\shell_1 - \shell_0$ as the discrete
counterpart to $\tau \dot \shell(0)$ for $\tau = \frac1K$.  Motivated by the
fact that $ \frac1K \; \log_{\shell_A}(\shell_B) = \frac1K \dot \shell(0)$ we
hence give the following definition of a discrete logarithm.  Suppose we have a
unique discrete geodesic $(\shell_0= \shell_A, \ldots, \shell_K= \shell_B)$,
then we define the \emph{discrete logarithm} by
\begin{align*}
  \discLog^{K}_{\shell_A}(\shell_B) = K (\shell_1 - \shell_0)\,.
\end{align*}
As in the continuous case, the discrete logarithm can be considered as the
linear representation of the nonlinear variation $\shell_B$ to $\shell_A$.

Next, we ask for fixed $\shell_0=\shell_A$ and
$\shell_1=\shell_A+\tfrac{\xi}{K}$ for a solution $\shell_2$ of the
Euler-Lagrange equation~\eqref{eq:EL} for $K=2$:
$0 = \energy_{,2}[\shell_0,\shell_1] + \energy_{,1}[\shell_1,\shell_2]\,.$ One
can show (cf.~\cite{RuWi12b}) that for $\shell_1 - \shell_0$ sufficiently small
a unique solution $\shell_2$ exists and thus $(\shell_0,\shell_1,\shell_2)$ is a
discrete geodesic. Iterating this procedure, we compute for given $\shell_{k-1}$
and $\shell_k$ a solution $\shell_{k+1}$ of the equation \eqref{eq:EL} and
obtain a variational scheme for the shooting of a discrete geodesics
$(\shell_0,\ldots, \shell_k,\ldots)$.  Thereby, $\shell_k$ is a discrete counter
part of $\exp_{\shell_A} \left(\tfrac{k}{K} \xi \right)$ with $\tau=\tfrac1K$
being the time step size of the scheme.  Hence, we introduce the \emph{discrete
  exponential map}
\begin{align*}
  \mathbf{Exp}^{K}_{\shell_A}(\tfrac{k}{K},\xi)= \shell_k\,.
\end{align*}

\begin{figure}
  \begin{center}
    \includegraphics[width=\linewidth]{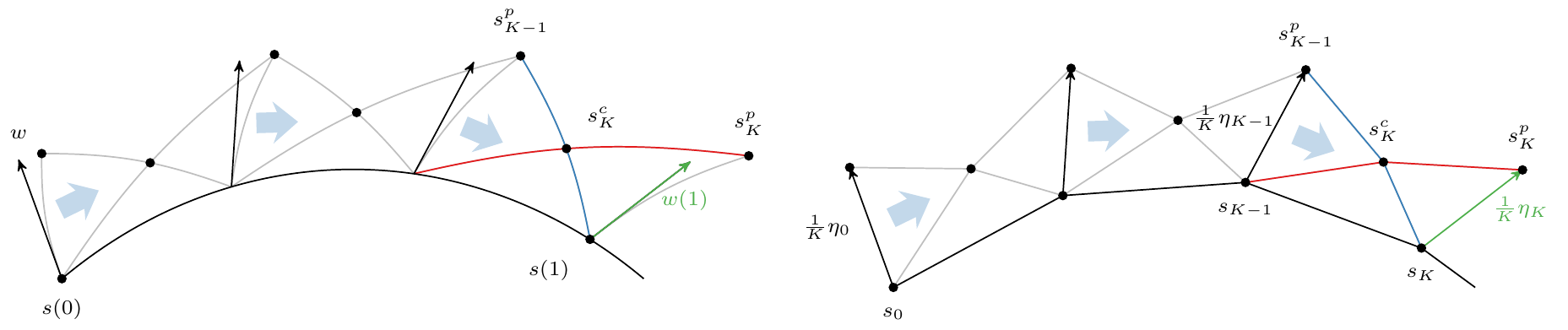}
  \end{center}
  \caption{A sketch of the discretization of parallel transport via Schild's
    ladder in the continuous context (left) and the discrete context (right).}
  \label{fig:conceptDiscreteParallelTransport}
\end{figure}

Let us emphasize that we use here a slightly different notation as to
\cite{RuWi12b} with a more natural scaling of arguments and results of the
discrete logarithm and the discrete exponential map.  Finally, let us recall the
definition of a discrete counterpart of the parallel transport.  \emph{Schild's
  ladder} (\cite{EhPiSc72,KhMiNe00}) is a well-known first-order approximation
of parallel transport along curves on Riemannian manifolds. It is based on the
construction of a sequence of geodesic parallelograms (see Figure
\ref{fig:conceptDiscreteParallelTransport}).  Given a path
$(\shell(t))_{t\in[0,1]}$ and a short tangent vector
$w=w_0 \in T_{\shell(0)}\shellspace$ the approximation of the parallel
transported vector at time $k\tau$ via a geodesic parallelogram can be expressed
by the iterative scheme
\begin{align*}
  \shell_{k-1}^p &= \exp_{\shell((k-1)\tau)}(w_{k-1}),\quad   
                   \shell_{k}^c = \exp_{\shell_{k-1}^p}(\tfrac{1}{2} \log_{\shell_{k-1}^p}(\shell(k\tau)) )\\
  \shell_{k}^p &= \exp_{\shell((k-1)\tau)}( 2 \log_{\shell((k-1)\tau)}(\shell_{k}^c) ), \text{  and } 
                 w_k = \log_{\shell(k\tau)}(\shell_{k}^p)
\end{align*}
for $k=1,\ldots, K$ and $w_K$ as an approximation of the parallel transport
$P_{(\shell(t))_{t\in[0,1]}} (w)$ (cf.\
Fig.~\ref{fig:conceptDiscreteParallelTransport}).  Now, we replace the
continuous operators $\exp$ and $\log$ by their discrete counterparts and obtain
the following iteration for a given discrete curve
$(\shell_0, \ldots, \shell_K)$ (with $\shell_k - \shell_{k-1}$ sufficiently
small for $k = 1, \ldots, K$) and a displacement $\eta = \eta_0$:
\begin{align*}
  \shell_{k-1}^p &= \shell_{k-1} + \tfrac{1}{K}\,\eta_{k-1} \\
  \shell_{k}^c &= \shell_{k-1}^p + \tfrac12 \discLog^{2}_{\shell_{k-1}^p}(\shell_k) \\ 
  \shell_{k}^p &= \discExp^{2}_{\shell_{k-1}}(1, \shell_{k}^c - \shell_{k-1} ) \\ 
  \eta_k &= K(\shell_{k}^p - \shell_{k}),
\end{align*}
for $k = 1, \ldots, K$ and define the discrete parallel transport of $\eta$
along $(\shell_0, \ldots, \shell_K)$ as (cf.\
Fig.~\ref{fig:conceptDiscreteParallelTransport} with marked objects for the case
$k=K$)
\begin{align*}
  \Ptransport_{(\shell_0,\ldots,\shell_K)}^K (\eta_0) = \eta_K.
\end{align*}
The scaling by $\tfrac1K$ in the first line and the rescaling by $K$ in the last
line is necessary to ensure consistency with the continuous parallel transport
for $K\to \infty$ and for $O(\eta_o) = 1$ as discussed in \cite{RuWi12b}.  In
contrast to the continuous setting, where $w_k \in T_{\shell_k}\shellspace$ are
tangent vectors, in the discrete setting $\eta_k = \shell_k^p - \shell_k$ are
displacements and $\Ptransport^K_{\shell_{k+1},\shell_k}(\eta_k) = \eta_{k+1}$
maps a displacement $\eta_k$ or the shell $\shell_{k}$ to a displacement
$\eta_{k+1}$ of the shell $\shell_{k+1}$.  The well-definedness of all these
operators and their first order convergence with respect to the time step size
$\tau=\tfrac1K$ is proved in~\cite{RuWi12b} under assumptions fulfilled in the
space of subdivision shell surfaces.

\begin{figure}[h]
  \centering    \includegraphics[width=1.0\linewidth]{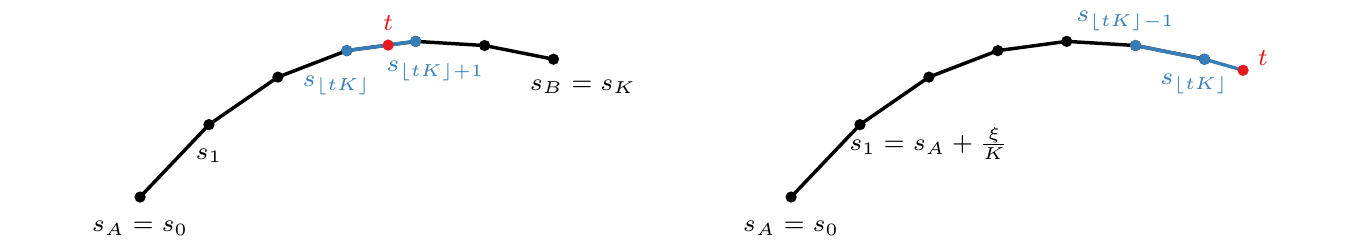}
  \caption{Sketch for the definition of the generalized interpolation operator
    (left) and the generalized exponential map (right). }
  \label{fig:general_operators}
\end{figure}

\paragraph{More general discrete interpolation and extrapolation}
In order to define subdivision curves we would like to compute discrete geodesic
averaging not only at discrete time steps \(\tfrac{k}{K}\) for
\(k = 0, \dots, K\) but also for any arbitrary time \(t \in [0,1]\) with $t$ not
being a multiple of $\tfrac1K$. Therefore, we generalize the discrete
interpolation operator in~\eqref{eq:1} as follows. We first compute a discrete
geodesic of length $K+1$ as explained above.  Then in a second step we use a
weighted, discrete geodesic interpolation between the two shells
$\shell_{\floor{tK}}$ and $\shell_{\floor{tK}+1}$ for which
$t\in (\tfrac{\floor{tK}}{K},\tfrac{\floor{tK}+1}{K})$ and $\floor{y}$ denotes
the largest integer smaller than or equal to $y \in \R$ and
$\shell_{k} = \Average^K\left(\shell_A, \shell_B, \tfrac{k}{K}\right)$.  Thereby
the weighted, discrete geodesic interpolation is given as
\begin{equation}
  \label{eq:2}
  \Average^K(\shell_A, \shell_B, t) =
  \argmin_{\shell \in \shellspace} (1-t(K)) \energy[\shell_{\floor{tK}}, \shell]
  + t(K) \energy[\shell, \shell_{\floor{tK}+1}], 
\end{equation}
where $t(K) = tK - \floor{tK}$ (cf.\ Fig.~\ref{fig:general_operators}); e.g.\
for $t=\frac{2k+1}{2K}$ and $t(K) = \frac12$ we retrieve the equal weighting
from the definition of a $3$ shell discrete geodesic.

In a similar manner we define a generalized exponential map for $t$ not being a
multiple of the time step size $\tau =\tfrac1K$ for given $K$.  To this end, we
define $\discExp^K(t,\xi)$ such that $\shell_k := \discExp^K(k/K,\xi)$ for
$k=\floor{tK}$ is the weighted discrete geodesic average of
$\shell_{k-1} = \discExp^K((k-1)/K,\xi)$ and the unknown shell
$\discExp^K(t,\xi)$ with weights $\tfrac{t(K)}{1+t(K)}$ and $\tfrac{1}{1+t(K)}$,
i.e.\
\begin{equation*}
  \shell_k = \argmin_{\tilde \shell} \left( \frac{t(K)}{1+t(K)} \energy[\shell_{k-1}, \tilde \shell]  + \frac{1}{1+t(K)}\energy[\tilde \shell,  \shell]\right)
\end{equation*}
with unknown shell $\shell = \discExp^K(t,\xi)$. This can be expressed in terms
of the Euler-Lagrange equation 
\begin{equation*}
  0= \frac{t(K)}{1+t(K)} \energy_{,2}[\shell_{k-1}, \shell_k]  + \frac{1}{1+t(K)} \energy_{,1}[\shell_k,  \shell]\,.  
\end{equation*}
The solution is unique for
sufficiently small $\xi$ (cf.\ Fig.~\ref{fig:general_operators}).  Obviously,
for $t(K) =0$ we obtain $\discExp^K(t,\xi) = \shell_k$ and for $t(K) =1$ and
thus $t=\tfrac{k+1}{K}$ the exponential map coincides with the previous
definition of $\discExp^K(\tfrac{k+1}{K},\xi)$.

\paragraph{Implementation}
For the discretization of shells in \(\shellspace\) we implemented a subdivision
element method in C++ based on Loop's subdivision functions. This scheme ensures
\(H^2\) regularity (\(C^2\) apart from irregular vertices) of the limit
function. Thus it allows for the energy densities used in the variational
approximation of the path energy a conforming discretization in particular of
the relative shape operator (cf.~\cite{Ci00}).  The crucial step in the
implementation of subdivision methods is the choice of the numerical
quadrature. Here, we use the mid-edge quadrature considered in~\cite{JuMa16}
which displays a good compromise between efficiency and
robustness. Additionally, this quadrature rule in conjunction with
lookup tables allows for the simulation with input meshes containing more than
one extraordinary vertex per patch. We refer to~\cite{JuMa16} for
implementational details.

\begin{figure}
  \begin{center}
    \includegraphics[width=\linewidth]{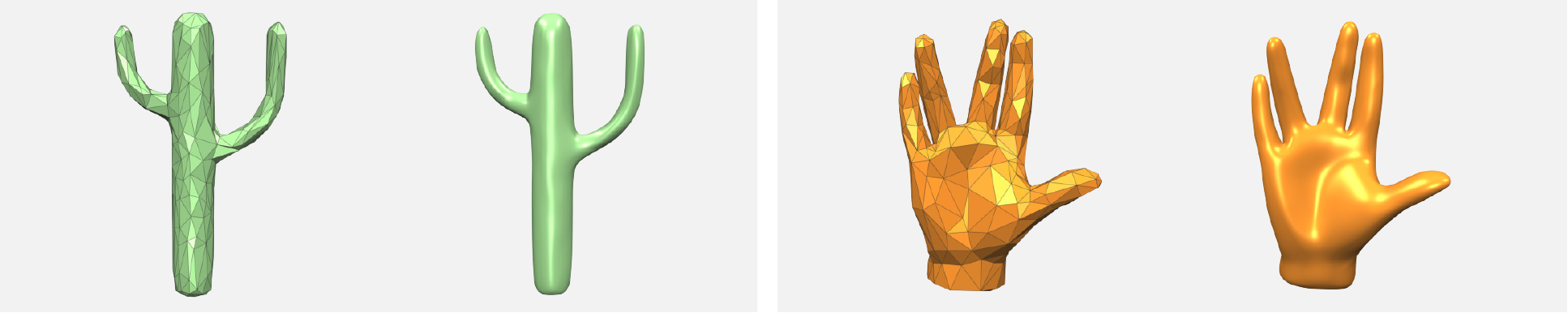}
  \end{center}
  \caption{Control mesh and subdivision limit surface for the cactus (left) and
    the hand shape (right).}
  \label{fig:Loop}
\end{figure}

For the computations presented in the subsequent sections we used a cactus and
a hand model (cf.\ Fig.~\ref{fig:Loop} for an instance of the control meshes and
the subdivision surfaces). The underlying meshes are closed and consist of 263
vertices for the cactus model and 305 vertices for the hand model.

The numerical computation of discrete geodesic, discrete logarithm, discrete
exponential map and discrete parallel transport uses the elastic deformation
energy defined in~\eqref{eq:totalEnergy} with membrane and bending energy
density given as in~\eqref{eq:membraneDensity}
and~\eqref{eq:bendingDensity}. In all computations the material parameters are
chosen as follows: \(\lambda = \mu = 1.0\) and \(\delta = 0.01\).
The nonlinear minimization of the discrete path energy~\eqref{eq:discPathEnergy}
for the computation of discrete geodesics and geodesic extrapolation is
performed by solving the corresponding set of Euler-Lagrange equations via
Newton's method with stepsize control. The iteration is stopped if the squared
\(\ell^2\)-norm of the Newton step decreases below \(\epsilon = 10^{-4}\).
Note that the performance of Newton's method strongly depends on the
initialization. While this does not pose a problem for the geodesic
extrapolation where the unknown shell \(\shell_k\) can be simply initialized
with \(\shell_{k-1}\) the issue is more delicate for geodesic interpolation. For
the computation of discrete geodesics between to given shells \(\shell_0\) and
\(\shell_K\) we therefore iteratively minimized the functionals 
\begin{equation*}
  \shell \mapsto (K-k) \energy[\shell_0, \shell] + k \energy[\shell, \shell_K]
\end{equation*}
for \(k = 1, \dots, K-1\) to get \(\shell_k\), starting the \(k\)th interation
with the inialization \(\shell = \shell_{k-1}\). This approach yields a good
initialization \((\shell_0, \dots, \shell_K)\) for the final Newton iteration
such that in most cases Newton's method is in the contraction region of
quadratic convergence.

Using this numerical setup the computation times needed in the simulations for
this paper range from a couple of minutes for the computation of simple quadratic
B\'ezier curves to several hours for subdivision curves of levels higher than
four.

\section{B\'ezier curves in shell space}
\label{sec:bezier-curves-shell}

Now, we consider B\'ezier curves in the space of subdivision shell surfaces and
derive a discrete counterpart of Riemannian B\'ezier curves using the
discretization from~\cite{HeRuWa12} and following~\cite{BrTyHi16} adapted to the
space of subdivision shell surfaces.  Consider a set of control shells
$(\shell_0^0, \ldots, \shell_n^0)$ with $\shell_j^0 \in \shellspace$ for
$j = 0,\ldots,n$ and the mapping
\begin{align*}
  \calB: \underbrace{\shellspace \times \ldots \times \shellspace}_{n+1} \times [0,1] \to \shellspace  
\end{align*}
which is recursively defined via the \emph{de Casteljau algorithm}
\begin{align*}
  \calB(\shell_i, \ldots, \shell_j, t) = \average( \calB(\shell_i, \ldots, \shell_{j-1}, t), \calB(\shell_{i+1}, \ldots, \shell_j, t), t  )
\end{align*}
for $i,j \in \{0,\ldots,n\}$ with $i < j$ and $\calB(\shell,t) = \shell$.  In
other words for $(\shell_0^0, \ldots, \shell_n^0)$ and fixed $t \in [0,1]$ we
compute for $j = 1, \ldots, n$ and $i = j, \ldots, n$ the shells
\begin{align*}
  \shell_{i}^{j} &= \average( \shell_{i-1}^{j-1}, \shell_{i}^{j-1}, t )
\end{align*}
and after $n$ steps we obtain the shell
$\calB(\shell_0^0, \ldots, \shell_n^0, t) = \shell_{n}^{n}$.  The resulting
curve $(\calB(\shell_0^0, \ldots, \shell_n^0, t))_{t \in [0,1]}$ is denoted a
\emph{B\'ezier curve in shell space of degree $n$}.  Compared to the Euclidean
case the shells $\shell_i^j$ do not lie on a straight line but on a geodesic
$( \average( \shell_{i-1}^{j-1}, \shell_{i}^{j-1}, t ) )_{t\in[0,1]}$ connecting
$\shell_{i-1}^{j-1}$ with $\shell_{i}^{j-1}$.

\begin{figure}
  \begin{center}
    \includegraphics[width=\linewidth]{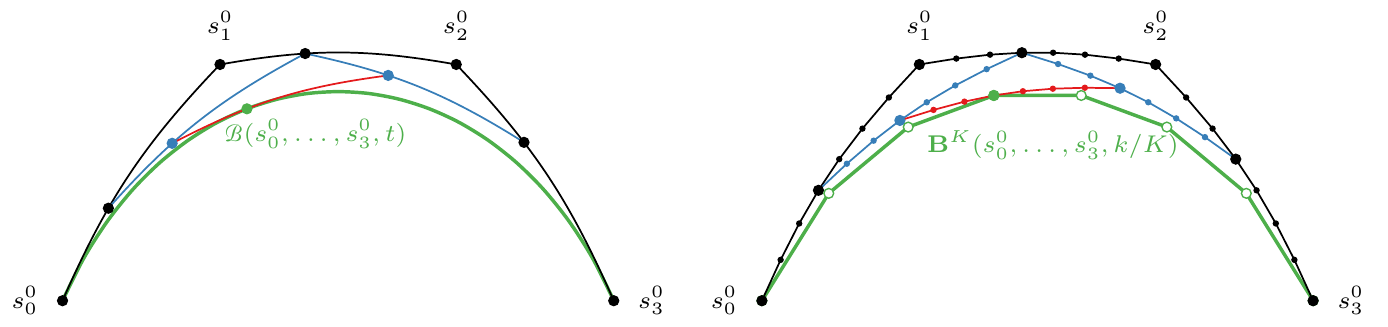}
  \end{center}
  \caption[]{Sketch of the continuous (left) and discrete (right) Riemannian de
    Casteljau algorithm defining a continuous and a discrete cubic B\'ezier
    curve, respectively (the algorithm proceeds as follows: \tikz{\node[circle,
      fill=black, inner sep=2pt] at (0,0) {}; \node[circle, fill=blue, inner
      sep=2pt] at (0.3,0) {}; \node[circle, fill=red, inner sep=2pt] at (0.6,0)
      {}; \node[circle, fill=green, inner sep=2pt] at (0.9,0) {};}).}
  \label{fig:conceptcubBezier}
\end{figure}

\begin{figure}[h]
  \centering
  \includegraphics[trim={0.8cm 0.6cm 0.8cm 0.8cm},
  width=\linewidth]{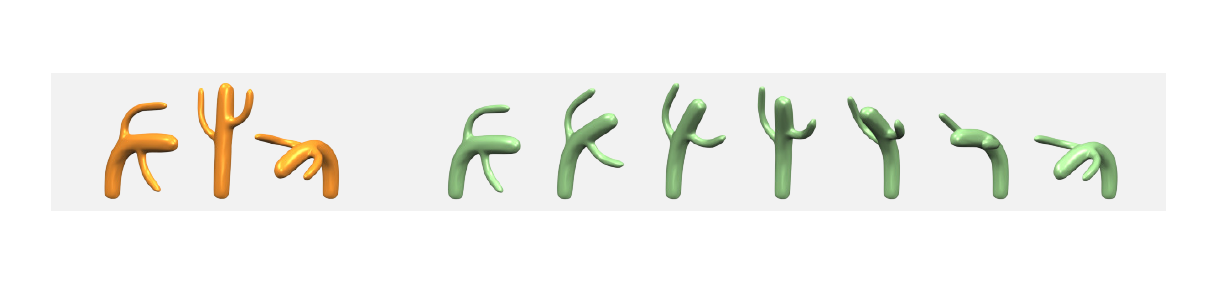}
  \includegraphics[trim={0.8cm 0.8cm 0.8cm 0.7cm},
  width=\linewidth]{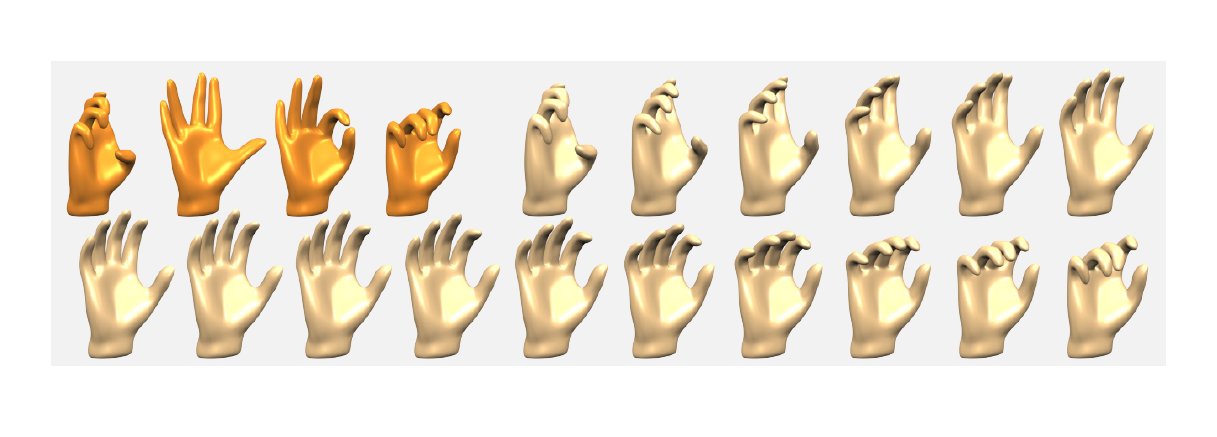}
  \caption{Two discrete B{\'e}zier curves for the cactus ($K=6$, green) and the
    hand model ($K=15$, beige).  The control shells are shown in orange on the
    left.}
  \label{fig:bezier_results}
\end{figure}

\begin{figure}[h]
  \centering
  \includegraphics[width=0.4\linewidth]{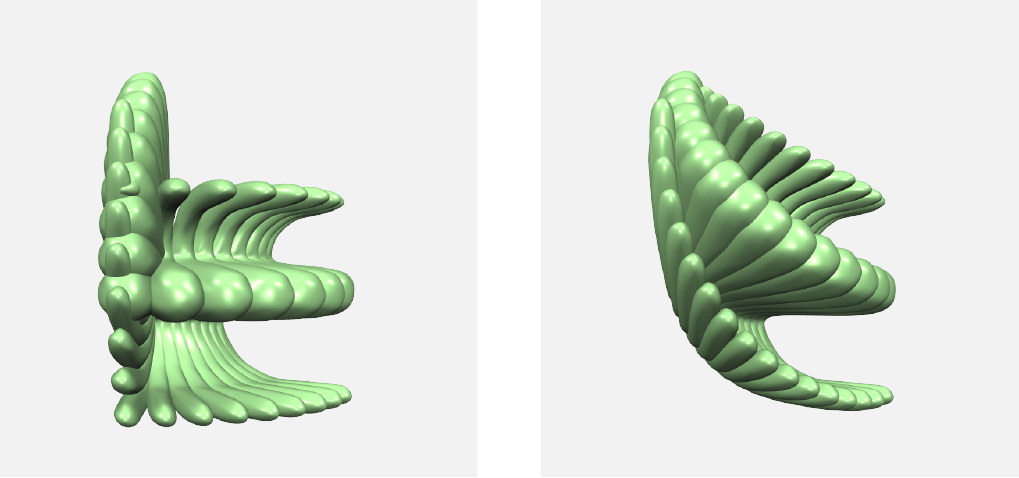}
  \caption{Comparison between a piecewise geodesic paths (left) and a quadratic
    B\'ezier curve (right).}
  \label{fig:bezier_comparison}
\end{figure}

Now, we transfer the concept of Riemannian B\'ezier curves to our time-discrete
setup.  Again, we consider a set of control shells
$(\shell_0^0, \ldots, \shell_n^0)$ with $\shell_j^0 \in \shellspace$ for
$j = 0,\ldots,n$ and some $K\in \N$.  Then we define a mapping
\begin{align*}
  \discB{K}: \underbrace{\shellspace \times \ldots \times \shellspace}_{n+1} \times \{0, \tfrac{1}{K}, \ldots, \tfrac{K}{K}\} \to \shellspace  
\end{align*}
recursively via the following \emph{discrete de Casteljau algorithm}
\begin{align*}
  \discB{K}(\shell_i, \ldots, \shell_j, \tfrac{k}{K}) = \Average^K( \discB{K}(\shell_i, \ldots, \shell_{j-1},  \tfrac{k}{K}), \discB{K}(\shell_{i+1}, \ldots, \shell_j, \tfrac{k}{K}),  \tfrac{k}{K} ) 
\end{align*}
for $i,j \in \{0,\ldots,n\}$ with $i < j$, $k \in \{0,\ldots,K\}$, and
$\discB{K}(\shell, \tfrac{k}{K}) = \shell$ .  For the evaluation of a B\'ezier
curve of degree $n$ with input shells $(\shell_0^0,\ldots,\shell_n^0)$, and for
$K\in \N$ the discrete de Casteljau algorithm for $k\in \{0,\ldots, K\}$ reads
as follows
\begin{center}
  \begin{algorithmic}
    \FOR{$j=1$ \TO $n$} \FOR{$i=j$ \TO $n$}
    \STATE{$\shell^j_i=
      \Average^K(\shell^{j-1}_{i-1},\shell^{j-1}_{i},\tfrac{k}{K})$}
    \ENDFOR
    \ENDFOR
    \STATE{$\discB{K}(\shell_0^0, \ldots,\shell_n^0,\tfrac{k}{K}) =
      \shell_n^n$}.
  \end{algorithmic}
\end{center}
In analogy to the continuous B\'ezier curve we call the resulting discrete path
\begin{align*}
  (\discB{K}(\shell_0, \ldots, \shell_n,  \tfrac{k}{K}))_{k \in \{0, \ldots, K\}}
\end{align*}
a \emph{discrete B\'ezier curve in shell space of degree $n$} (see
Figures~\ref{fig:conceptcubBezier},\ref{fig:bezier_results},
and~\ref{fig:bezier_comparison}).  Let us remark that the evaluation of the
B\'ezier curve for general times $t\in [0,1]$, which are not multiples of
$\tfrac1K$ is straightforward using the the generalized interpolation.

\section{Time-discrete cardinal splines in shell space}
\label{sec:time-discr-card}

\begin{figure}
  \begin{center}
    \includegraphics[width=\linewidth]{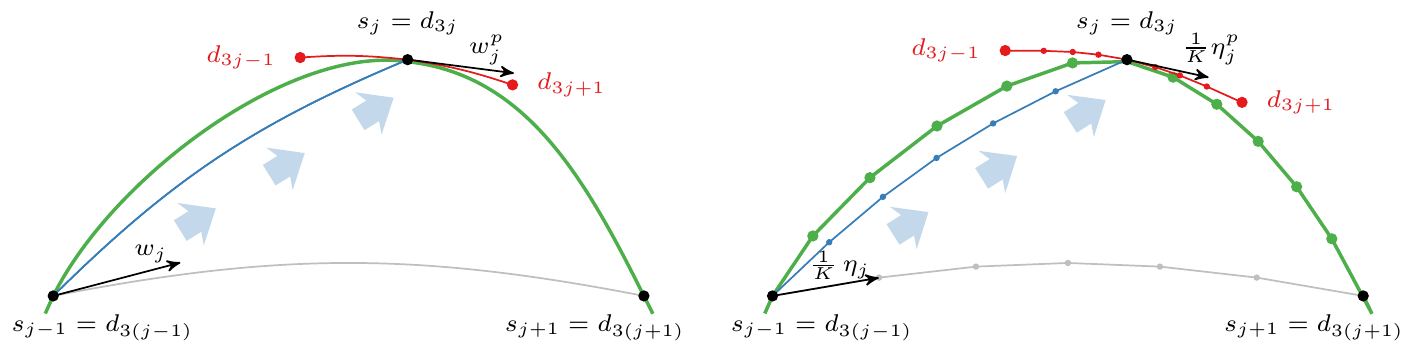}
  \end{center}
  \caption[]{Sketch of the construction of the control shells $\control_{3j-1}$
    and $\control_{3j+1}$ for the cardinal spline in the continuous (left) and
    discrete (right) setting (the algorithm proceeds as follows:
    \tikz{\node[circle, fill=black, inner sep=2pt] at (0,0) {}; \node[circle,
      fill=blue!30, inner sep=2pt] at (0.3,0) {}; \node[circle, fill=red, inner
      sep=2pt] at (0.6,0) {}; \node[circle, fill=green, inner sep=2pt] at
      (0.9,0) {};}).}
  \label{fig:conceptDiscreteCardinalSpline}
\end{figure}

In what follows we will use the relation between cubic B\'ezier curves and cubic
Hermite curves~\cite{Fa02} to define discrete cardinal splines in this space.
To this end, we first recall that a cubic Hermite curve
$t \mapsto \calH(\shell_A, v_A, v_B, \shell_B, t)$ with endpoints $\shell_A$ and
$\shell_B$ and with initial and final velocity $v_A\in T_{\shell_A}\shellspace$
and $v_B \in T_{\shell_B}\shellspace$, respectively, coincides with the
B{\'e}zier curve $t \mapsto \calB(\shell_0, \shell_1, \shell_2, \shell_3,t)$,
where
\begin{align*}
  \shell_0 &= \shell_A\,, \quad 
             \shell_1 = \exp_{\shell_A}(\tfrac{1}{3} v_A)\,, \quad
             \shell_2 = \exp_{\shell_B}(-\tfrac{1}{3} v_B)\,, \quad 
             \shell_3 = \shell_B\,.
\end{align*}
From this we immediately derive the discrete counterpart and define the
\emph{discrete cubic Hermite curve}
$$(\discH{K}(\shell_A, \xi_A, \xi_B, \shell_B, \tfrac{k}{K}))_{k \in \{0,\ldots,K\}}$$ for
given shells $\shell_A, \shell_B \in \shellspace$, displacements $\xi_A$,
$\xi_B$ as approximations of the tangent vectors $v_A$ and $v_B$, respectively,
and $K > 1$ by
\begin{align*}
  \discH{K}(\shell_A, \xi_A, \xi_B, \shell_B, \tfrac{k}{K} ) = \discB{K}(\shell_0, \shell_1, \shell_2, \shell_3, \tfrac{k}{K} ) \quad \text{where} \\ 
  \shell_0 = \shell_A\,, \quad  
  \shell_1 = \discExp^K_{\shell_A}(\frac{1}{3},\xi_A)\,, \quad
  \shell_2 = \discExp^K_{\shell_B}(\frac{1}{3},-\xi_B)\,, \quad 
  \shell_3 = \shell_B.
\end{align*}
Let us remark that for $k$ being a multiple of $3$ the term
$\discExp^K_{\shell_A}(\frac{1}{3},\xi_A)$ requires $\tfrac{K}3$ steps of the
iterative scheme for the discrete exponential map. Now, we are in the position
to discuss continuous cardinal splines on the Riemannian manifold $\shellspace$
of shells and to derive a proper \emph{discrete cardinal spline}. Cardinal
splines interpolate the given sequence of control shells.  In between each
consecutive pair of control points the cardinal spline coincides with a cubic
Hermite curve with the same pair of end shells.  The velocities at the end
shells of such a Hermite curve are computed via parallel transport of the initial
velocity of the geodesic pointing from the previous to the next shell.  An
exception is only the velocity at the first and the last shell.  The latter
mechanism ensures differentiability of the resulting spline curve. In detail,
for a sequence of shells $(\shell_0, \ldots, \shell_m)$ with
$\shell_j \in \shellspace$ for $j = 0,\ldots,m$, $m \geq 2$ and
$\kappa \in [0,3]$ consider the mapping
\begin{align*}
  \CS_{\kappa}: \underbrace{\shellspace \times \ldots \times \shellspace}_{m+1} \times [0,m] \to \shellspace  
\end{align*}
which is defined piecewise for $t \in [l, l+1]$ where $l = 0, \ldots, m-1$ by
\begin{align*}
  \CS_{\kappa}( \shell_0, \ldots, \shell_m, t ) = \calB( \control_{3l}, \control_{3l+1}, \control_{3l+2}, \control_{3l+3}, t-l )\,.
\end{align*}
To ensure the interpolation condition we have to set
\begin{align*}
  \control_{3j} &= \shell_j, \qquad j = 0, \ldots, m\,.
\end{align*}
On the left and on the right we prescribe boundary conditions via
\begin{align*}
  \control_{1} &= \average( \shell_0, \shell_1, \rfrac{\kappa}{3} ), \\
  \control_{3m-1} &= \average( \shell_{m-1}, \shell_m, 1 - \rfrac{\kappa}{3} ) = \average( \shell_m, \shell_{m-1}, \rfrac{\kappa}{3} )\,. 
\end{align*}
In explicit, this ensures that
$\frac{\mathrm{d}}{\mathrm{d} t} \CS_{\kappa}( \shell_0, \ldots, \shell_m, t )$
is parallel to $\log_{d_0} (d_1)$ for $t=0$ and to $\log_{d_{3m}} (d_{3m-1})$
for $t=m$\,.  Finally, we prescribe for each $j\in 1,\ldots, m-1$ the control
points $\control_{3j-1}$ and $\control_{3j+1}$\,.  To this end, we first compute
the initial velocity \(w_j\) of the geodesic from $\control_{3(j-1)}$ and
$\control_{3(j+1)}$ and scale it with $\rfrac{\kappa}{3}$, where $\kappa$ is the
tension factor, which controls the inflection behaviour of the cardinal spline:
\begin{align*}
  w_j &= \tfrac{\kappa}{3} \log_{\shell_{j-1}}(\shell_{j+1})\,. 
\end{align*} 
Then we transport this velocity along the geodesic from $\shell_{j-1}$ to
$\shell_{j}$ and use the exponential map in the positive and negative direction
of the transported velocity field to define the remaining control points:
\begin{align*}
  w_j^p &= \ptransport_{(\average(\shell_{j-1},\shell_{j},t))_{t \in [0,1]}} ( w_j ), \\
  \control_{3j-1} &= \exp_{\shell_j}(- w_j^p), \\
  \control_{3j+1} &= \exp_{\shell_j}( w_j^p)\,.
\end{align*}
For the computation of the shells $\control_{3j-1}$ and $\control_{3j+1}$ we
refer to the sketch in Figure~\ref{fig:conceptDiscreteCardinalSpline}.  The
resulting curve
$$(\CS_{\kappa}( \shell_0, \ldots, \shell_m, t ))_{t\in[0,m]}$$
defines the \emph{cubic cardinal spline in shell space} with \emph{tension
  parameter} $\kappa \in [0,3]$\,.

\begin{figure}[h]
  \centering
  \includegraphics[trim={0.8cm 0.6cm 0.8cm 0.8cm},
  width=\linewidth]{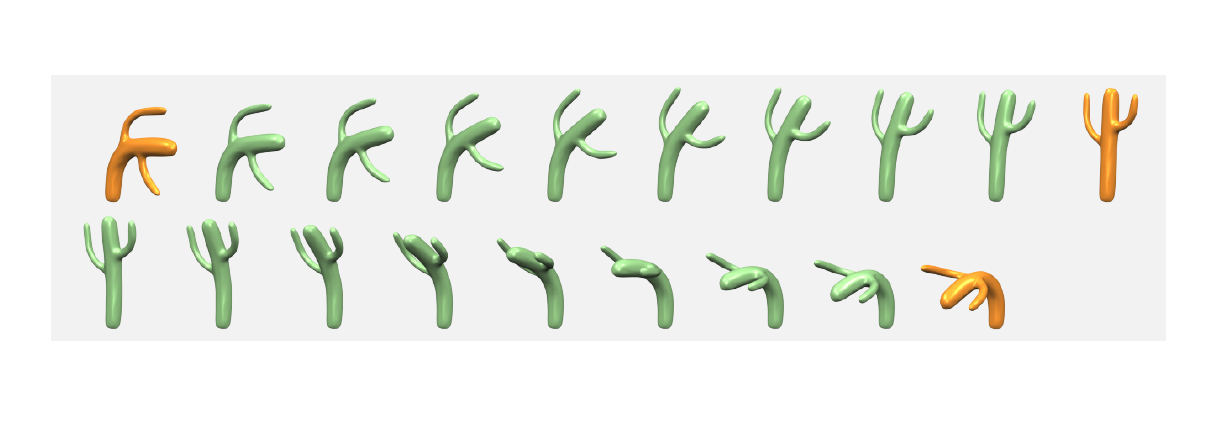}
  \includegraphics[trim={0.8cm 0.6cm 0.8cm 0.8cm},
  width=\linewidth]{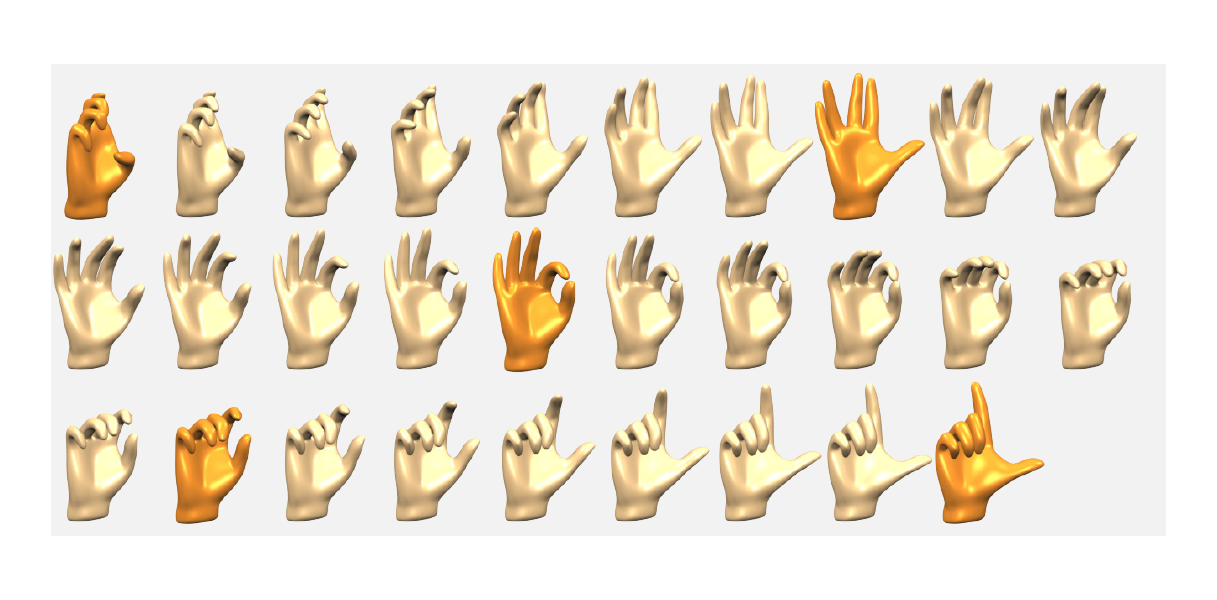}
  \caption{Two discrete cardinal spline curves for the cactus (green, \(K=9\)) and the
    hand model (beige, \(K=7\)).  The interpolated key frames are shown in
    orange.}
  \label{fig:cardinal_results}
\end{figure}

This construction can be transferred to the discrete context to define the
corresponding discrete cardinal spline.  Again we take into account a sequence
of control shells $(\shell_0, \ldots, \shell_m)$ with $\shell_j \in \shellspace$
for $j = 0,\ldots,m$, $m \geq 2$ and $K \geq 2$ and define the mapping
\begin{align*}
  \discCS{K}: \underbrace{\shellspace \times \ldots \times \shellspace}_{m+1} \times \{0, \tfrac{1}{K}, \ldots, m \cdot \tfrac{K}{K}\} \to \shellspace  
\end{align*}
which is given piecewise for $k \in [ l \cdot K, \ldots, (l+1) \cdot K]$ as
$l = 0, \ldots, m-1$ by
\begin{align*}
  \discCS{K}( \shell_0, \ldots, \shell_m,\tfrac{k}{K} ) = \discB{K}( \control_{3l}, \control_{3l+1}, \control_{3l+2}, \control_{3l+3}, \tfrac{k}{K}-l)\,,
\end{align*}
where the interpolation control points are copied as before, i.e.
\begin{align*}
  \control_{3j} &= \shell_j, \qquad j = 0, \ldots, m\,.
\end{align*}
The discrete counterparts of the boundary conditions are
\begin{align*}
  \control_{1}    &= \discExp^K_{\shell_0}\left( \rfrac{\kappa}{3} ,\discLog^K_{\shell_0}(\shell_1) \right), \\
  \control_{3m-1} &=  \discExp^K_{\shell_m}\left( \rfrac{\kappa}{3}, \discLog^K_{\shell_m}(\shell_{m-1}) \right)\,. \\
\end{align*}
The control points $\control_{3j-1}$ and $\control_{3j+1}$ for
$j\in 1,\ldots, m-1$ are determined as follows:
\begin{align*}
  \eta_j &=   \discLog^K_{\shell_{j-1}}( \shell_{j+1} )\,,\\
  \eta_j^p &= \Ptransport^{K}_{(\Average^K(\shell_{j-1},\shell_{j},\tfrac{k}{K}))_{k = \{0, \ldots, K\}}} \left( \eta_j \right)\,, \\
  \control_{3j-1} &= \discExp^{K}_{\shell_j} \left( \tfrac{\kappa}{3}, -\eta_j^p \right)\,,\\
  \control_{3j+1} &= \discExp^{K}_{\shell_j}\left( \tfrac{\kappa}{3}, \eta_j^p \right)\,. \\
\end{align*}
(see also Figure~\ref{fig:conceptDiscreteCardinalSpline}).  We call the
resulting curve
$(\discCS{K}( \shell_0, \ldots, \shell_m, \tfrac{k}{K} ))_{k\in\{0,\ldots,m
  \cdot K\}}$ the \emph{discrete cubic cardinal spline in shell space} with
tension $\kappa \in [0,3]$.  Let us emphasize that this is equivalent to the
construction with discrete cubic Hermite splines where a segment consists of the
shells $\shell_j$ and $\shell_{j+1}$ and displacements $\eta_{j-1}^p$ and
$-\eta_{j}^p$. Examples are given in Figure~\ref{fig:cardinal_results}.

\section{Interpolatory subdivision curves in shell space}
\label{sec:interp-subd-curv}

In this section we transfer the concept of unitary subdivision to shell space in
order to compute smooth interpolation curves for a given set of key frame
shells.  In Euclidean space a subdivision curve is defined as the limit of
increasingly fine control polygons subject to a small set of subdivision
rules. Replacing straight lines by geodesic paths, this construction can be
transferred to the Riemannian setting. For this purpose we make use of the
extension of the geodesic averaging operator and the exponential map in
Section~\ref{sec:vari-time-discr}.  Furthermore, for the ease of presentation
we generalize the averaging $\average$ to an interpolation $\interpolation$,
which is also defined for \(t \not\in [0,1]\).  In explicit, $\interpolation$ is
given as follows:
\begin{alignat*}{2}
  \interpolation(\shell_A, \shell_B, t) &= \average(\shell_A, \shell_B, t)
  &&\qquad\quad \mbox{ for } 0\leq t \leq 1\,,\\
  \interpolation(\shell_A, \shell_B, t) &= \exp_{\shell_A} (t
  \log_{\shell_A}(\shell_B))
  &&\qquad\quad \mbox{ for } t <0\,,\\
  \interpolation(\shell_A, \shell_B, t) &= \exp_{\shell_B} (-(t-1)
  \log_{\shell_B}(\shell_A)) &&\qquad\quad \mbox{ for } t >1\,.
\end{alignat*}
The same notational generalization can be performed for the discrete averaging:
\begin{alignat*}{2}
  \Interpolation^K(\shell_A, \shell_B, t) &= \Average^K(\shell_A, \shell_B, t)
  &&\quad \mbox{ for } 0\leq t \leq 1\,,\\
  \Interpolation^K(\shell_A, \shell_B, t) &= \discExp^K_{\shell_A} (-t,
  -\discLog^K_{\shell_A} (\shell_B))
  &&\quad \mbox{ for } t <0\,,\\
  \Interpolation^K(\shell_A, \shell_B, t) &= \discExp^K_{\shell_B} ((t-1), -
  \discLog^K_{\shell_B} (\shell_A)) &&\quad \mbox{ for } t >1\,.
\end{alignat*}
In the following we state the subdivision rules for the interpolatory binary
four- and six-point schemes as well as the interpolatory ternary four-point
scheme. Let us emphasize that these rules depends on a proper combination of
(discrete) geodesic interpolation and extrapolation via the (discrete)
exponential map, both encoded in the generalized interpolation operators
$\interpolation$ and $\Interpolation^K$.

\begin{figure}[h]
  \centering
  \includegraphics[width=\linewidth]{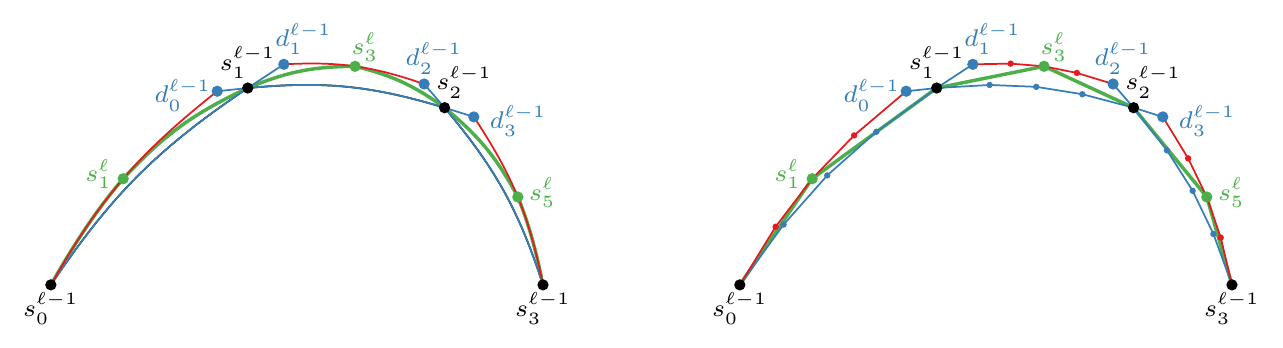}
  \caption[]{Sketch of the construction for the interpolatory binary four-point
    scheme in the continuous case (left) and the discrete case (right) (the
    algorithm proceeds as follows: \tikz{\node[circle, fill=black, inner
      sep=2pt] at (0,0) {}; \node[circle, fill=blue, inner sep=2pt] at (0.3,0)
      {}; \node[circle, fill=red, inner sep=2pt] at (0.6,0) {}; \node[circle,
      fill=green, inner sep=2pt] at (0.9,0) {};}).}
  \label{fig:bin4}
\end{figure}

\begin{figure}[h]
  \centering
  \includegraphics[trim={0.8cm 0.8cm 0.8cm 0.8cm},
  width=\linewidth]{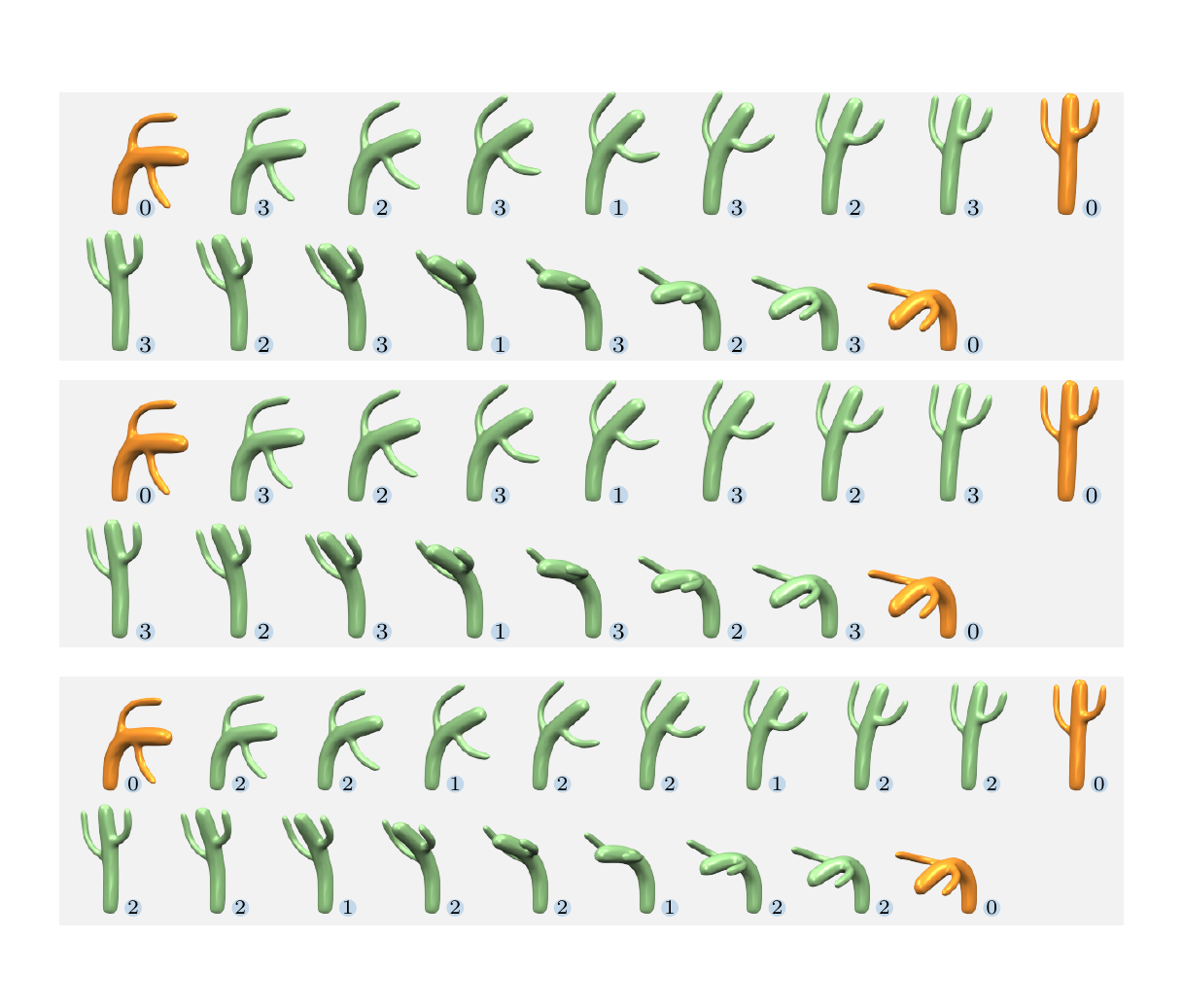}
  \caption{Different discrete subdivision curves interpolate three different
    poses of the cactus shell (orange): discrete binary four-point scheme (top),
    discrete binary six-point scheme (middle), discrete ternary four-point
    scheme.  The numbers indicate the level of refinement.}
  \label{fig:subdiv_cactus}
\end{figure}

\begin{figure}[h]
  \centering
  \includegraphics[width=0.8\linewidth]{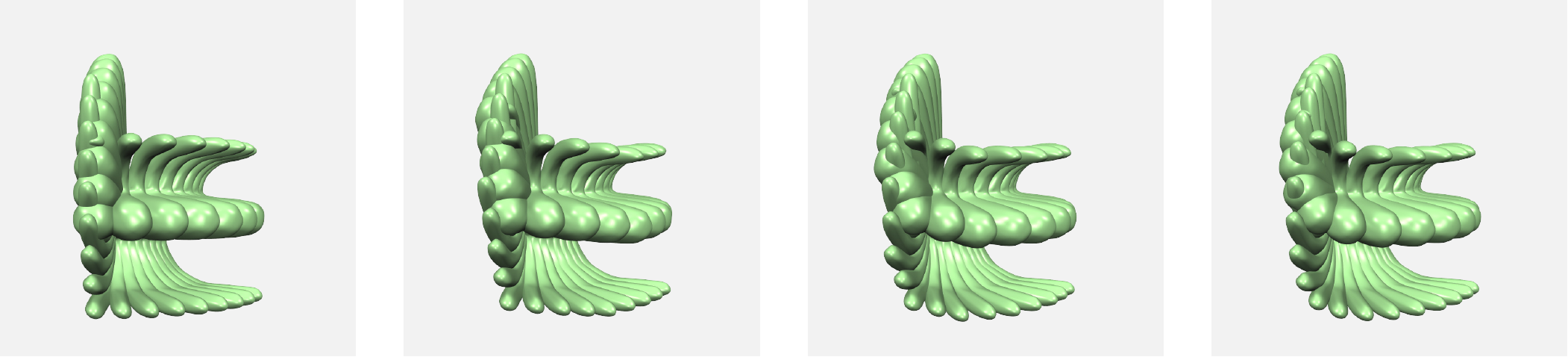}
  \caption{Comparison for binary 4-point rule for level
    $l= 0,\, 1, \,2, \,3$\,.}
  \label{fig:bin4_comparison}
\end{figure}

\begin{figure}[h]
  \centering
  \includegraphics[trim={0.8cm 0.8cm 0.8cm 0.8cm},
  width=\linewidth]{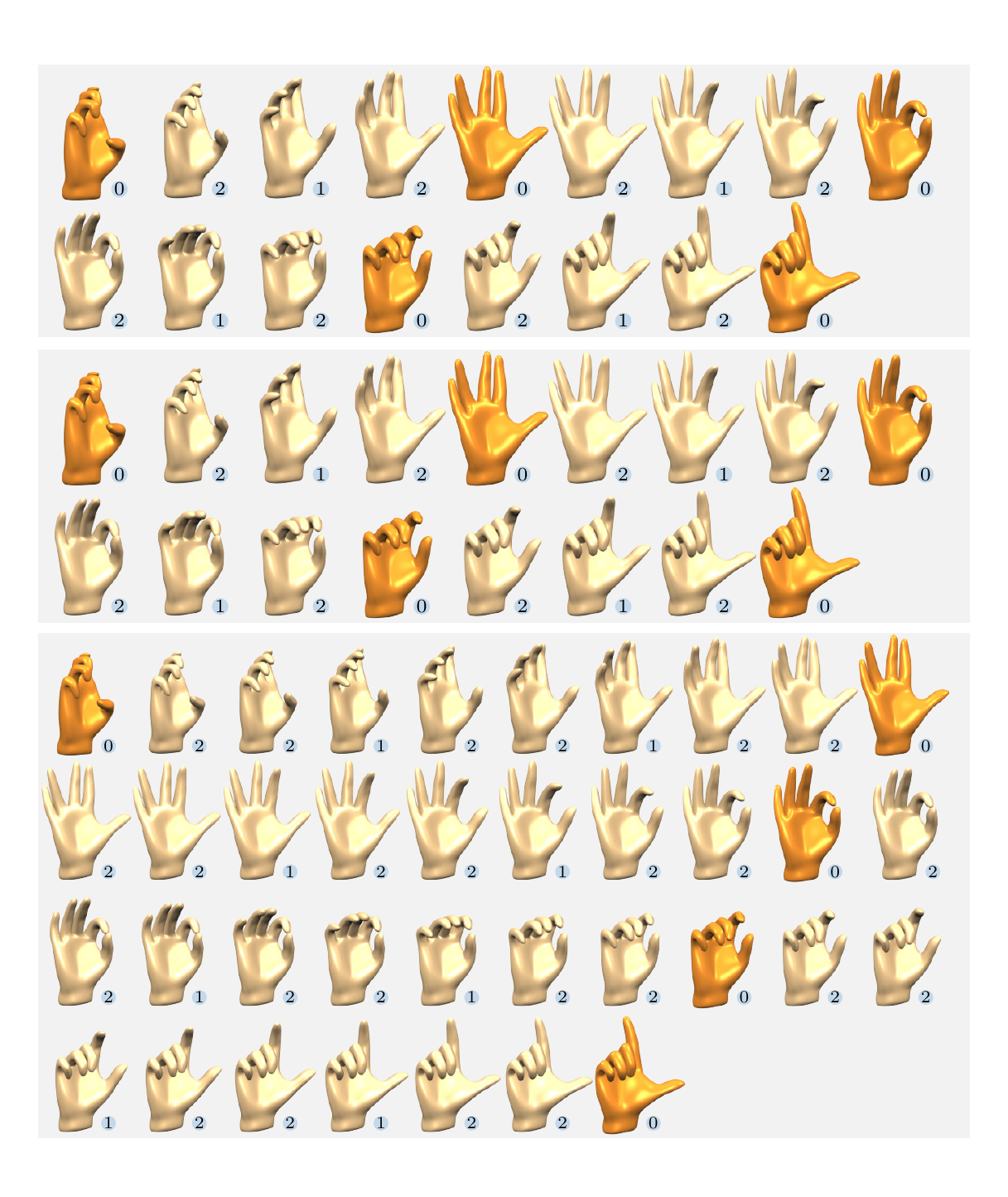}
  \caption{As above in Fig.~\ref{fig:subdiv_cactus} different discrete
    subdivision curves now interpolate five different poses of the hand model
    (orange): discrete binary four-point scheme (top), discrete binary six-point
    scheme (middle), discrete ternary four-point scheme.  Again the numbers
    indicate the level of refinement.}
  \label{fig:subdiv_hand}
\end{figure}

To this end, let us consider a set of initial shells
\(\shell_0^0, \dots, \shell_n^0 \in \shellspace\) with \(n \geq 2\).  The
\emph{Riemannian interpolatory binary four-point scheme} defines the
corresponding control shells \(\shell_k^\ell\) of level \(\ell > 0\) iteratively
by
\begin{align*}
  \shell_{2k}^\ell &= \shell_k^{\ell-1}, \quad 
                     \shell_{2k+1}^\ell = \interpolation(d_{2k}^{\ell-1}, d_{2k+1}^{\ell-1}, \tfrac{1}{2}), 
\end{align*}
for \(k = 0, \dots, 2^\ell n\) with
\begin{align*}
  d_{2k}^{\ell-1} &= \interpolation(\shell_{k-1}^{\ell-1}, \shell_k^{\ell-1}, \tfrac{9}{8}),  \quad 
                    d_{2k+1}^{\ell-1} = \interpolation(\shell_{k+2}^{\ell-1}, \shell_{k+1}^{\ell-1}, \tfrac{9}{8})
\end{align*}
(cf.\ the sketch in Fig.~\ref{fig:bin4}).  Similarly the \emph{Riemannian
  interpolatory binary six-point scheme} is given by
\begin{align*}
  \shell_{2k}^\ell &= \shell_k^{\ell-1},  \quad 
                     \shell_{2k+1}^\ell = \interpolation(d_{4k+2}^{\ell-1}, d_{4k+3}^{\ell-1},\tfrac{1}{2}),
\end{align*}
for \(k = 0, \dots, 2^\ell n\), where
\begin{alignat*}{2}
  d_{4k}^{\ell-1} &= \interpolation(\shell_{k-2}^{\ell-1}, \shell_{k-1}^{\ell-1}, \tfrac{25}{22}), \quad  
  d_{4k+1}^{\ell-1} &&=  \interpolation(\shell_{k+3}^{\ell-1}, \shell_{k+2}^{\ell-1}, \tfrac{25}{22}), \\
  d_{4k+2}^{\ell-1} &= \interpolation(d_{4k}^{\ell-1}, \shell_k^{\ell-1}, \tfrac{75}{64}), \quad 
  d_{4k+3}^{\ell-1} &&= \interpolation(d_{4k+1}^{\ell-1} \shell_{k+1}^{\ell-1}, \tfrac{75}{64}).
\end{alignat*}
Finally the \emph{interpolatory ternary four-point scheme} consists of the rules
\begin{align*}
  \shell_{3k}^\ell &= \shell_k^{\ell-1}, \quad
  \shell_{3k+1}^\ell = \interpolation(d_{4k}^{\ell-1}, d_{4k+1}^{\ell-1}, \tfrac{10}{33}), \quad
  \shell_{3k+2}^\ell = \interpolation(d_{4k+2}^{\ell-1}, d_{4k+3}^{\ell-1}, \tfrac{10}{33}), 
\end{align*}
with \(k = 0, \dots, 3^\ell n\), and
\begin{alignat*}{2}
  d_{4k}^{\ell-1} &= \interpolation(\shell_{k-1}^{\ell-1}, \shell_k^{\ell-1}, \tfrac{76}{69}), \quad
                    d_{4k+1}^{\ell-1} &&= \interpolation(\shell_{k+1}^{\ell-1}, \shell_{k+2}^{\ell-1}, -\tfrac{2}{15}), \\
  d_{4k+2}^{\ell-1} &= \interpolation(\shell_{k+2}^{\ell-1}, \shell_{k+1}^{\ell-1}, \tfrac{76}{69}), \quad
                      d_{4k+3}^{\ell-1} &&= \interpolation(\shell_k^{\ell-1}, \shell_{k-1}^{\ell-1}, -\tfrac{2}{15}). 
\end{alignat*}
In order to transfer these subdivision rules into the time-discrete setup we
replace the time-continuous (generalized) interpolation operator by its
time-discrete counterpart.  For \(K \geq 2\) and for the set of keyframe shells
\(\shell_0^0, \dots, \shell_n^0 \in \shellspace\) (our control shells) with
\(n \geq 2\) we obtain the \emph{interpolatory, discrete Riemannian binary
  four-point scheme} with shells \(\shell_k^\ell\) on level \(\ell > 0\)
iteratively defined by
\begin{align*}
  \shell_{2k}^\ell &= \shell_k^{\ell-1}, \quad 
                     \shell_{2k+1}^\ell = \Interpolation^K(d_{2k}^{\ell-1}, d_{2k+1}^{\ell-1}, \tfrac{1}{2}), 
\end{align*}
for \(k = 0, \dots, 2^\ell n\) with
\begin{align*}
  d_{2k}^{\ell-1} &= \Interpolation^K(\shell_{k-1}^{\ell-1}, \shell_k^{\ell-1}, \tfrac{9}{8}),  \quad 
                    d_{2k+1}^{\ell-1} = \Interpolation^K(\shell_{k+2}^{\ell-1}, \shell_{k+1}^{\ell-1}, \tfrac{9}{8})
\end{align*}
(cf.\ the sketch in Fig.~\ref{fig:bin4}).  The \emph{discrete binary six-point
  scheme} reads as follows
\begin{align*}
  \shell_{2k}^\ell &= \shell_k^{\ell-1},  \quad 
                     \shell_{2k+1}^\ell = \Interpolation^K(d_{4k+2}^{\ell-1}, d_{4k+3}^{\ell-1},\tfrac{1}{2}),
\end{align*}
for \(k = 0, \dots, 2^\ell n\), where
\begin{alignat*}{2}
  d_{4k}^{\ell-1} &= \Interpolation^K(\shell_{k-2}^{\ell-1}, \shell_{k-1}^{\ell-1}, \tfrac{25}{22}), \quad  
  d_{4k+1}^{\ell-1} &&=  \Interpolation^K(\shell_{k+3}^{\ell-1}, \shell_{k+2}^{\ell-1}, \tfrac{25}{22}),\\
  d_{4k+2}^{\ell-1} &= \Interpolation^K(d_{4k}^{\ell-1}, \shell_k^{\ell-1}, \tfrac{75}{64}), \quad 
  d_{4k+3}^{\ell-1} &&= \Interpolation^K(d_{4k+1}^{\ell-1} \shell_{k+1}^{\ell-1}, \tfrac{75}{64}).
\end{alignat*}
For the \emph{discrete ternary four-point scheme} we get
\begin{align*}
  \shell_{3k}^\ell &= \shell_k^{\ell-1}, \quad
  \shell_{3k+1}^\ell = \Interpolation^K(d_{4k}^{\ell-1}, d_{4k+1}^{\ell-1}, \tfrac{10}{33}), \quad
  \shell_{3k+2}^\ell = \Interpolation^K(d_{4k+2}^{\ell-1}, d_{4k+3}^{\ell-1}, \tfrac{10}{33}), 
\end{align*}
with \(k = 0, \dots, 3^\ell n\), and
\begin{alignat*}{2}
  d_{4k}^{\ell-1} &= \Interpolation^K(\shell_{k-1}^{\ell-1}, \shell_k^{\ell-1}, \tfrac{76}{69}), \quad
  d_{4k+1}^{\ell-1} && = \Interpolation^K(\shell_{k+1}^{\ell-1}, \shell_{k+2}^{\ell-1}, -\tfrac{2}{15}), \\
  d_{4k+2}^{\ell-1} &= \Interpolation^K(\shell_{k+2}^{\ell-1}, \shell_{k+1}^{\ell-1}, \tfrac{76}{69}), \quad
  d_{4k+3}^{\ell-1} &&= \Interpolation^K(\shell_k^{\ell-1}, \shell_{k-1}^{\ell-1}, -\tfrac{2}{15}). 
\end{alignat*}
Examples for the different subdivision schemes in shell space are given in
Figures~\ref{fig:subdiv_cactus} and~\ref{fig:subdiv_hand} and a comparison is
shown in Figure~\ref{fig:bin4_comparison}.

\section{Discussion}
\label{sec:discussion}

We have generalized the notion of different interpolating curves, namely
cardinal splines, and three different types of interpolatory subdivision schemes
for curves to the Riemannian manifold of shells. This methodology turns into a
computationally efficient toolbox via the implementation of a variational time
discretization of geodesic interpolation, extrapolation via a discrete geometric
exponential map, discrete geometric logarithm, and discrete parallel transport.
Due to the physical foundation of the underlying Riemannian metric on the space
of shells the application to subdivision shell surfaces led to physically
intuitive interpolation curves in shell space.

The use of subdivision surfaces to represent shells allows a conforming
discretization of the path energy and does not require any discrete geometry
version of the relative shape operator.  
In particular, approximating a given smooth surface with subdivision surfaces generated by an 
increasing number of control points we can expect convergence of the relative shape operator.
Concerning the numerical implementation
Newton's method is used in particular to solve the geodesic interpolation
problem, which requires the minimization of a nonlinear energy over $K-1$
coupled control meshes of the subdivision shells.  Here the hierarchical
approach incorporated in the subdivision method is particularly well-suited. It
ensures that Newton's method is in most cases (especially on finer levels of the
subdivision scheme) in the contraction region of quadratic convergence.

We remark that the construction of the interpolating curves presented in
this paper is entirely based on the solution of problems that are local in
time. However, the corresponding functionals subject to minimization are
non-convex such that the Newton solver may terminate in local minima which
could possibly lead to discontinuities and non-uniqueness of the resulting
interpolation curves. In the application this problem can be prevented by a good
initialization of the Newton iteration (like the iterative initialization
presented in Section~\ref{sec:vari-time-discr}) and a suitable choice of the
step size \(\tau = 1/K\). We never encountered the described problem in any of
our computations.

Compared to the Riemannian spline approach by Heeren et al.~\cite{HeRuSc16} we
have to solve solely local problems in time, whereas in~\cite{HeRuSc16} a fully
coupled global optimization problem with a large set of nonlinear constraints
has to be solved.  Furthermore, we do not require the simplification via the
$L\Theta A$ approach for more complex shell models, which requires a back
projection onto the manifold of shell surfaces.  On the other hand our approach
follows the paradigms of constructive approximation and lacks an energy
minimizing interpretation.

Furthermore, the convergence behaviour of the sequence of subdivision curves to
a smooth limit surfaces is open in particular for fixed time step size
$\tau =\tfrac1K$. The proximity condition approach by Wallner and Dyn
~\cite{WaDy05} might be a good starting points to study smoothness
properties. In our case two types of perturbations have to be considered: the
perturbation due to the curviness of the manifold and the perturbation caused by
our variational approximation of geodesics.  Recently, Wallner~\cite{Wa14}
proved convergence of the linear four-point scheme and other univariate
interpolatory schemes in Riemannian manifolds.  To this end he studied
Riemannian edge length contractivity of the schemes on a manifold combined with
a multiresolutional analysis.  Indeed, to apply these results certain
consistency results of the path length of geodesic edges and the corresponding
discrete path length are already known (cf.~\cite{RuWi12b}).

\section*{Acknowledgement}
We acknowledge support by the FWF in Austria under the grant S117 (NFN) and the
DFG in Germany under the grant Ru 567/14-1.

\bibliographystyle{acm}
\bibliography{bibtex/own,bibtex/all,bibtex/library,newBibTexItems}

\def\polhk#1{\setbox0=\hbox{#1}{\ooalign{\hidewidth
  \lower1.5ex\hbox{`}\hidewidth\crcr\unhbox0}}}
\begin{thebibliography}{10}

\bibitem{AbGoSt16}
{\sc Absil, P.-A., Gousenbourger, P.-Y., Striewski, P., and Wirth, B.}
\newblock Differentiable piecewise-{B}\'{e}zier surfaces on {R}iemannian
  manifolds.
\newblock {\em SIAM Journal on Imaging Sciences\/} (2016).
\newblock accepted.

\bibitem{AdOvWa08}
{\sc Adams, B., Ovsjanikov, M., Wand, M., Seidel, H.-P., and Guibas, L.~J.}
\newblock Meshless modeling of deformable shapes and their motion.
\newblock {\em Proceedings of the 2008 ACM SIGGRAPH/Eurographics Symposium on
  Computer Animation\/} (2008).

\bibitem{Ar01}
{\sc Arden, G.}
\newblock {\em Approximation Properties of Subdivision Surfaces}.
\newblock PhD thesis, University of Washington, 2001.

\bibitem{BaBr11}
{\sc Bauer, M., and Bruveris, M.}
\newblock A new {R}iemannian setting for surface registration.
\newblock In {\em Proc. of MICCAI Workshop on Mathematical Foundations of
  Computational Anatomy\/} (2011), pp.~182--194.
\newblock arXiv:1106.0620.

\bibitem{BaHaMi11}
{\sc Bauer, M., Harms, P., and Michor, P.~W.}
\newblock Sobolev metrics on shape space of surfaces.
\newblock {\em J. Geom. Mech. 3}, 4 (2011), 389--438.

\bibitem{BrTyHi16}
{\sc Brandt, C., von Tycowicz, C., and Hildebrandt, K.}
\newblock Geometric flows of curves in shape space for processing motion of
  deformable objects.
\newblock {\em Comput. Graph. Forum 35}, 2 (2016), 295--305.

\bibitem{BrBrKi09}
{\sc Bronstein, A.~M., Bronstein, M.~M., and Kimmel, R.}
\newblock Topology-invariant similarity of nonrigid shapes.
\newblock {\em Inter. J. Comput. Vision 81}, 3 (2009), 281--301.

\bibitem{CaSiCr01}
{\sc Camarinha, M., Silva~Leite, F., and Crouch, P.}
\newblock On the geometry of {R}iemannian cubic polynomials.
\newblock {\em Diff. Geom. Appl. 15\/} (2001), 107--135.

\bibitem{Ca12}
{\sc Cashman, T.~J.}
\newblock {B}eyond {C}atmull-{C}lark? {A} survey of advances in subdivision
  surface methods.
\newblock {\em Computer Graphics Forum 31}, 1 (2012), 42--61.

\bibitem{CaCl78}
{\sc Catmull, E., and Clark, J.}
\newblock {R}ecursively generated {B}-spline surfaces on arbitrary topological
  meshes.
\newblock {\em Computer Aided Design 10\/} (1978), 350--355.

\bibitem{Ci00}
{\sc Ciarlet, P.~G.}
\newblock {\em Mathematical elasticity. {V}ol. {III}}, vol.~29 of {\em Studies
  in Mathematics and its Applications}.
\newblock North-Holland Publishing Co., Amsterdam, 2000.
\newblock Theory of shells.

\bibitem{CiLo11}
{\sc Cirak, F., and Long, Q.}
\newblock Subdivision shells with exact boundary control and non--manifold
  geometry.
\newblock {\em Internat. J. Numer. Methods Engrg. 88}, 9 (2011), 897--923.

\bibitem{CiOr01}
{\sc Cirak, F., and Ortiz, M.}
\newblock {F}ully {C}1--conforming subdivision elements for finite deformation
  thin--shell analysis.
\newblock {\em Internat. J. Numer. Methods Engrg. 51(7)\/} (2001), 813--833.

\bibitem{CiOrSc00}
{\sc Cirak, F., Ortiz, M., and Schr{\"{o}}der, P.}
\newblock Subdivision surfaces: a new paradigm for thin-shell finite--element
  analysis.
\newblock {\em Internat. J. Numer. Methods Engrg. 47}, 12 (2000), 2039--72.

\bibitem{CiScAn02}
{\sc Cirak, F., Scott, M.~J., Antonsson, E.~K., Ortiz, M., and Schr{\"{o}}der,
  P.}
\newblock Integrated modeling, finite--element analysis and engineering design
  for thin--shell structures using subdivision.
\newblock {\em Comput.-Aided Des. 34}, 2 (2002), 137--148.

\bibitem{CrSi95}
{\sc Crouch, P., and Silva~Leite, F.}
\newblock The dynamic interpolation problem: On {R}iemannian manifold, {L}ie
  groups and symmetric spaces.
\newblock {\em J. Dynam. Control Systems 1}, 2 (1995), 177--202.

\bibitem{DeKaTr98}
{\sc DeRose, T., Kass, M., and Truong, T.}
\newblock {S}ubdivision {S}urfaces in {C}haracter {A}nimation.
\newblock In {\em Proceedings of the 25th Annual Conference on Computer
  Graphics and Interactive Techniques\/} (1998), SIGGRAPH '98, pp.~85--94.

\bibitem{DoCarmo1992}
{\sc do~Carmo, M.~P.}
\newblock {\em Riemannian geometry}.
\newblock Birkh{\"{a}}user Boston Inc., 1992.
\newblock Translated from the second Portuguese edition.

\bibitem{DoSa78}
{\sc Doo, D., and Sabin, M.}
\newblock Behaviour of recursive division surfaces near extraordinary points.
\newblock {\em Computer-Aided Design 10(6)\/} (1978), 356--360.

\bibitem{DuXi13}
{\sc Duchamp, T., Xie, G., and Yu, T.}
\newblock Single basepoint subdivision schemes for manifold-valued data:
  time-symmetry without space-symmetry.
\newblock {\em Found. Comput. Math. 13}, 5 (2013), 693--728.

\bibitem{Dy02}
{\sc Dyn, N.}
\newblock Interpolatory subdivision schemes.
\newblock In {\em Tutorials on Multiresolution in Geometric Modelling},
  A.~Iske, E.~Quak, and M.~S. Floater, Eds. Springer, Berlin, 2002, pp.~25--50.

\bibitem{DyGr10}
{\sc Dyn, N., Grohs, P., and Wallner, J.}
\newblock Approximation order of interpolatory nonlinear subdivision schemes.
\newblock {\em J. Comput. Appl. Math. 233}, 7 (2010), 1697--1703.

\bibitem{EfRuSi14}
{\sc Effland, A., Rumpf, M., Simon, S., Stahn, K., and Wirth, B.}
\newblock B\'{e}zier curves in the space of images.
\newblock In {\em Proc. of International Conference on Scale Space and
  Variational Methods in Computer Vision}, vol.~9087 of {\em Lecture Notes in
  Computer Science}. Springer, Cham, 2015, pp.~372--384.

\bibitem{EhPiSc72}
{\sc Ehlers, J., Pirani, F. A.~E., and Schild, A.}
\newblock The geometry of free fall and light propagation.
\newblock In {\em General relativity (papers in honour of {J}. {L}. {S}ynge)}.
  Clarendon Press, Oxford, 1972, pp.~63--84.

\bibitem{Fa02}
{\sc Farin, G.}
\newblock {\em Curves and Surfaces for CAGD: A Practical Guide}, 5th~ed.
\newblock Morgan Kaufmann Publishers Inc., San Francisco, CA, USA, 2002.

\bibitem{FrBo11}
{\sc Fr{\"{o}}hlich, S., and Botsch, M.}
\newblock Example-driven deformations based on discrete shells.
\newblock {\em Comput. Graph. Forum 30}, 8 (2011), 2246--2257.

\bibitem{GiGi02}
{\sc Giambo, R., and Giannoni, F.}
\newblock An analytical theory for {R}iemmanian cubic polynomials.
\newblock {\em IMA Journal of Mathematical Control and Information 19\/}
  (2002), 445--460.

\bibitem{GoSaAb14}
{\sc Gousenbourger, P.-Y., Samir, C., and Absil, P.}
\newblock Piecewise-bezier c1 interpolation on riemannian manifolds with
  application to 2d shape morphing.
\newblock In {\em Pattern Recognition (ICPR), 2014 22nd International
  Conference on\/} (2014), pp.~4086--4091.

\bibitem{GrTu04}
{\sc Green, S., and Turkiyyah, G.}
\newblock {S}econd {O}rder {A}ccurate {C}onstraint {F}ormulation for
  {S}ubdivision {F}inite {E}lement {S}imulation of {T}hin {S}hells.
\newblock {\em International Journal For Numerical Methods In Engineering 61},
  3 (2004), 380--405.

\bibitem{GrTuSt02}
{\sc Green, S., Turkiyyah, G., and Storti, D.}
\newblock {S}ubdivision--{B}ased {M}ultilevel {M}ethods for {L}arge {S}cale
  {E}ngineering {S}imulation of {T}hin {S}hells.
\newblock In {\em Proceedings of ACM Solid Modeling\/} (2002).

\bibitem{Gr10}
{\sc Grohs, P.}
\newblock A general proximity analysis of nonlinear subdivision schemes.
\newblock {\em SIAM J. Math. Anal. 42}, 2 (2010), 729--750.

\bibitem{HeRuSc14}
{\sc Heeren, B., Rumpf, M., Schr{\"{o}}der, P., Wardetzky, M., and Wirth, B.}
\newblock Exploring the geometry of the space of shells.
\newblock {\em Comput. Graph. Forum 33}, 5 (2014), 247--256.

\bibitem{HeRuSc16}
{\sc Heeren, B., Rumpf, M., Schr{\"{o}}der, P., Wardetzky, M., and Wirth, B.}
\newblock Splines in the space of shells.
\newblock {\em Comput. Graph. Forum 35}, 5 (2016), 111--120.

\bibitem{HeRuWa12}
{\sc Heeren, B., Rumpf, M., Wardetzky, M., and Wirth, B.}
\newblock Time-discrete geodesics in the space of shells.
\newblock {\em Comput. Graph. Forum 31}, 5 (2012), 1755--1764.

\bibitem{HiMuFl12}
{\sc Hinkle, J., Muralidharan, P., Fletcher, P.~T., and Joshi, S.}
\newblock Polynomial regression on {R}iemannian manifolds.
\newblock In {\em Proc. of European Conference on Computer Vision\/} (2012),
  vol.~7574 of {\em Lecture Notes in Computer Science}, pp.~1--14.

\bibitem{HoPo04}
{\sc Hofer, M., and Pottmann, H.}
\newblock Energy-minimizing splines in manifolds.
\newblock {\em ACM Trans. Graph. 23}, 3 (2004), 284--293.

\bibitem{JiZeLu09}
{\sc Jin, M., Zeng, W., Luo, F., and Gu, X.}
\newblock Computing teichm\"uller shape space.
\newblock {\em IEEE Transactions on Visualization and Computer Graphics 15}, 3
  (2009), 504--517.

\bibitem{JuMa16}
{\sc J{\"u}ttler, B., Mantzaflaris, A., Perl, R., and Rumpf, M.}
\newblock On numerical integration in isogeometric subdivision methods for
  {PDEs} on surfaces.
\newblock {\em Computer Methods in Applied Mechanics and Engineering 302\/}
  (2016), 131--146.

\bibitem{KaAn08}
{\sc Kass, M., and Anderson, J.}
\newblock Animating oscillatory motion with overlap: Wiggly splines.
\newblock {\em ACM Trans. Graph. 27}, 3 (2008), 28:1--28:8.

\bibitem{KhMiNe00}
{\sc Kheyfets, A., Miller, W.~A., and Newton, G.~A.}
\newblock Schild's ladder parallel transport procedure for an arbitrary
  connection.
\newblock {\em Internat. J. Theoret. Phys. 39}, 12 (2000), 2891--2898.

\bibitem{KiMiPo07}
{\sc Kilian, M., Mitra, N.~J., and Pottmann, H.}
\newblock Geometric modeling in shape space.
\newblock {\em ACM Trans. Graph. 26\/} (2007), 1--8.

\bibitem{KuKlDi10}
{\sc Kurtek, S., Klassen, E., Ding, Z., and Srivastava, A.}
\newblock A novel {R}iemannian framework for shape analysis of {3D} objects.
\newblock In {\em Proc. of IEEE Conference on Computer Vision and Pattern
  Recognition\/} (2010), pp.~1625--1632.

\bibitem{KuKlGo12}
{\sc Kurtek, S., Klassen, E., Gore, J.~C., Ding, Z., and Srivastava, A.}
\newblock Elastic geodesic paths in shape space of parametrized surfaces.
\newblock {\em IEEE Trans. Pattern Anal. Mach. Intell. 34}, 9 (2012),
  1717--1730.

\bibitem{LiShDi10}
{\sc Liu, X., Shi, Y., Dinov, I., and Mio, W.}
\newblock A computational model of multidimensional shape.
\newblock {\em Int. J. Comput. Vis. 89}, 1 (2010), 69--83.

\bibitem{Lo94}
{\sc Loop, C.}
\newblock Smooth spline surfaces over irregular meshes.
\newblock In {\em SIGGRAPH '94: Proceedings of the 21st annual conference on
  Computer graphics and interactive techniques\/} (1994), ACM Press,
  pp.~303--310.

\bibitem{NoHePa89}
{\sc Noakes, L., Heinzinger, G., and Paden, B.}
\newblock Cubic splines on curved spaces.
\newblock {\em IMA J. Math. Control Inform. 6}, 4 (1989), 465--473.

\bibitem{PeRe08}
{\sc Peters, J., and Reif, U.}
\newblock {\em Subdivision Surfaces}.
\newblock Springer Series in Geometry and Computing, 2008.

\bibitem{PeWu06a}
{\sc Peters, J., and Wu, X.}
\newblock On the local linear independence of generalized subdivision
  functions.
\newblock {\em SIAM Journal on Numerical Analysis (SINUM) 44}, 6 (2006),
  2389--2407.

\bibitem{PoNo07}
{\sc Popiel, T., and Noakes, L.}
\newblock B{\'{e}}zier curves and ${C}^2$ interpolation in {R}iemannian
  manifolds.
\newblock {\em J. Approx. Theory 148}, 2 (2007), 111--127.

\bibitem{PoHo05}
{\sc Pottmann, H., and Hofer, M.}
\newblock A variational approach to spline curves on surface.
\newblock {\em Comput. Aided Geom. Design 22}, 7 (2005), 693--709.

\bibitem{Re95}
{\sc Reif, U.}
\newblock A unified approach to subdivision algorithms near extraordinary
  vertices.
\newblock {\em Comput. Aided Geom. Design 12}, 2 (1995), 153--174.

\bibitem{ReSc01}
{\sc Reif, U., and Schr{\"{o}}der, P.}
\newblock Curvature integrability of subdivision surfaces.
\newblock {\em Adv. Comput. Math. 14}, 2 (2001), 157--174.

\bibitem{RuWi12b}
{\sc Rumpf, M., and Wirth, B.}
\newblock Variational time discretization of geodesic calculus.
\newblock {\em IMA J. Numer. Anal. 35}, 3 (2015), 1011--1046.

\bibitem{ScTySe14}
{\sc Schulz, C., von Tycowicz, C., Seidel, H.-P., and Hildebrandt, K.}
\newblock Animating deformable objects using sparse spacetime constraints.
\newblock {\em ACM Trans. Graph. 33}, 4 (July 2014), 109:1--109:10.

\bibitem{ScTySe15}
{\sc Schulz, C., von Tycowicz, C., Seidel, H.-P., and Hildebrandt, K.}
\newblock Animating articulated characters using wiggly splines.
\newblock In {\em Proceedings of the ACM SIGGRAPH / Eurographics Symposium on
  Computer Animation\/} (2015), pp.~101--109.

\bibitem{St99}
{\sc Stam, J.}
\newblock Evaluation of {L}oop subdivision surfaces.
\newblock SIGGRAPH 99 Course Notes, 1999.

\bibitem{ThWaSt06}
{\sc Thomaszewski, B., Wacker, M., and Stra{\ss}er, W.}
\newblock A consistent bending model for cloth simulation with corotational
  subdivision finite elements.
\newblock In {\em Proc. of ACM SIGGRAPH/Eurographics Symposium on Computer
  Animation\/} (2006), pp.~107--116.

\bibitem{TrVi12}
{\sc Trouv{\'e}, A., and Vialard, F.-X.}
\newblock Shape splines and stochastic shape evolutions : A second order point
  of view.
\newblock {\em Quart. Appl. Math. 70}, 2 (2012), 219--251.

\bibitem{Wa04}
{\sc Wallner, J.}
\newblock Existence of set-interpolating and energy-minimizing curves.
\newblock {\em Comput. Aided Geom. Design 21}, 9 (2004), 883--892.

\bibitem{Wa06b}
{\sc Wallner, J.}
\newblock Smoothness analysis of subdivision schemes by proximity.
\newblock {\em Constr. Approx. 24}, 3 (2006), 289--318.

\bibitem{Wa14}
{\sc Wallner, J.}
\newblock On convergent interpolatory subdivision schemes in {R}iemannian
  geometry.
\newblock {\em Constr. Approx. 40}, 3 (2014), 473--486.

\bibitem{WaDy05}
{\sc Wallner, J., and Dyn, N.}
\newblock Convergence and {$C^1$} analysis of subdivision schemes on manifolds
  by proximity.
\newblock {\em Comput. Aided Geom. Design 22}, 7 (2005), 593--622.

\bibitem{WiDrAl10}
{\sc Winkler, T., Drieseberg, J., Alexa, M., and Hormann, K.}
\newblock Multi-scale geometry interpolation.
\newblock {\em Comput. Graph. Forum 29}, 2 (2010), 309 -- 318.
\newblock Proceedings of Eurographics.

\bibitem{WiKa88}
{\sc Witkin, A., and Kass, M.}
\newblock Spacetime constraints.
\newblock {\em Computer Graphics 22}, 4 (1988), 159--168.

\bibitem{ZoJKo14}
{\sc Zore, U., J{\"{u}}ttler, B., and Kosinka, J.}
\newblock {O}n the {L}inear {I}ndependence of ({T}runcated) {H}ierarchical
  {S}ubdivision {S}plines.
\newblock Tech. rep., NFN G+S: Technical Report No. 17, 2014.

\end{thebibliography}

\end{document}